\documentclass[a4paper]{article}

\title{Interaction patterns and individual dynamics shape the way we move in synchrony}
\date{}

\author{
Francesco Alderisio\footnotemark[2] , Gianfranco Fiore\footnotemark[2] , Robin N. Salesse\footnotemark[3] , \\ Beno{\^\i}t  G. Bardy\footnotemark[3] \footnotemark[4] \ \& Mario di Bernardo\footnotemark[2] \footnotemark[5] \ \text{*}
}

\usepackage{amsmath}
\usepackage{amssymb}
\usepackage{amsthm}
\usepackage{hyperref}
\usepackage{graphicx}
\usepackage{subfig}
\usepackage{epstopdf}
\usepackage{booktabs}
\usepackage{tabularx}  
\usepackage{xfrac}
\usepackage{color}
\usepackage{lineno}
\usepackage[margin=0.8in]{geometry}

\newtheorem{remark}{Remark}

\begin{document}

\begin{titlepage}

\maketitle
An important open problem in Human Behaviour is to understand how coordination emerges in human ensembles. This problem has been seldom studied quantitatively in the existing literature, in contrast to situations involving dual interaction. Here we study motor coordination (or synchronisation) in a group of individuals where participants are asked to visually coordinate an oscillatory hand motion. We separately tested two groups of seven participants. We observed that the coordination level of the ensemble depends on group homogeneity, as well as on the pattern of visual couplings (who looked at whom). Despite the complexity of social interactions, we show that networks of coupled heterogeneous oscillators with different structures capture well the group dynamics. Our findings are relevant to any activity requiring the coordination of several people, as in music, sport or at work, and can be extended to account for other perceptual forms of interaction such as sound or feel.

\footnotetext[2]{Department of Engineering Mathematics, Merchant Venturers Building, University of Bristol, Woodland Road, Clifton, Bristol BS8 1UB, United Kingdom (\texttt{f.alderisio@bristol.ac.uk}, \texttt{gianfranco.fiore@bristol.ac.uk}, \texttt{m.dibernardo@bristol.ac.uk})}
\footnotetext[3]{EuroMov, Montpellier University, 700 Avenue du Pic Saint-Loup, 34090 Montpellier, France (\texttt{benoit.bardy@umontpellier.fr}, \texttt{salesse.robin@gmail.com})}
\footnotetext[4]{Institut Universitaire de France, 1 rue Descartes, 75231 Paris Cedex 05, France}
\footnotetext[5]{Department of Electrical Engineering and Information Technology, University of Naples Federico II, Via Claudio 21, 80125 Naples, Italy (\texttt{mario.dibernardo@unina.it})}

\end{titlepage}

\setcounter{page}{2}

\section*{Introduction}
Motor coordination and synchronisation is an essential feature of many human activities, where a group of individuals performs a joint task. Examples include hands clapping in an audience \cite{NRBVB00}, walking in a crowd \cite{MPGHT10,RW14}, music playing \cite{DBLTCCAF12,BDGCF14}, sports \cite{DWA10,S16} or dance \cite{L14,DLTL16,EBQBM16}.
Achieving synchronisation in the group involves perceptual interaction through sound, feel, or sight, and the establishment of mental connectedness and social attachment among group members \cite{VRHKV04,WH09}.
This human phenomenon has rarely been studied in the existing literature, in contrast to the large number of results on the dynamics of animal groups \cite{CKFL05,NABV10,NVPMVB13,ZBPdB15}.

Indeed, most available theoretical results on human coordination are mostly concerned with the case of two individuals performing a joint action \cite{OdGJLK08,ST94,VMLB11,slowinski2015dynamic}, a recent example being that of the \emph{mirror game} \cite{NDA11}, presented as a paradigmatic case for the study of how people imitate each other's movements in a pair \cite{ZATdMSMC,ZASTdB16}.

For larger groups of individuals, available results are mostly experimental observations of group behaviour and include studies on rocking chairs \cite{FR10,RGFGM12,ABdB15}, rhythmic activities and marching tasks \cite{IR15}, choir singers during a concert \cite{HT09}, group synchronisation of arm movements and respiratory rhythms \cite{CBVB14}, team rowing during a race \cite{WW95} and a few other sport situations \cite{YY11}. These studies have analysed the emergent level of coordination in the group, but never in relation to the structure of interactions or the individual dynamics of group members. 

Other studies have shown that the outcome and the quality of the performance in a number of situations strongly depend on how the individuals in the ensemble exchange visual, auditory and motor information \cite{DBLTCCAF12,healey2005inter,kauffeld2009complaint,passos2011networks,duarte2012intra,RTAR13}. Here too, these studies lack information about how specific interaction patterns affect coordination in the group, and in general a systematic and quantitative evaluation is missing of how coupling structure and intrinsic homogeneity (or heterogeneity) in the group contribute to the emergence of synchronisation.

In this work, we address this open problem and \textbf{confirm for the first time, experimentally and computationally, that different visual interaction patterns in the group affect the coordination level achieved by its members.} We take as a paradigmatic example the case where participants are asked to generate an oscillatory hand motion and coordinate it with that of the others. In addition, \textbf{we unfold the effects on group synchronisation of heterogeneities in the individual motion characteristics of the participants} (measured in terms of the intrinsic frequency of oscillation they generate in isolation).

Specifically, we show that the level of coordination achieved by group members is influenced by the combined action of the features characterising their motion in isolation (i.e., their natural oscillation frequency) and the specific interconnections (i.e., topological structure) among the players. We find that some topologies (e.g., all-to-all) give rise to higher levels of synchronisation (defined as an overall reduction in the phase mismatch among individuals) regardless of individual differences, whereas for other topologies (e.g., consecutive dyads) a better synchronisation is achieved through a higher homogeneity in individual dynamics.

We also propose a data-driven mathematical model that captures most of the coordination features observed experimentally. The model shows that, surprisingly, when performing a simple oscillatory movement, the group behaves as a network of nonlinearly coupled heterogeneous oscillators \cite{K84,strogatz2000kuramoto} despite the complexity of unavoidable social interactions in the group \cite{BL95,MFH10,SFV13,VDBCF16}.
Also, the model reproduces the dependence of the coordination level of each individual in the group upon the intrinsic properties of its members and the interaction structure among them, notwithstanding the complex neural mechanisms behind the emergence of such coordination.

\section*{Results}
\begin{figure}[!ht]
 \centering
 \includegraphics[width=.75\textwidth]{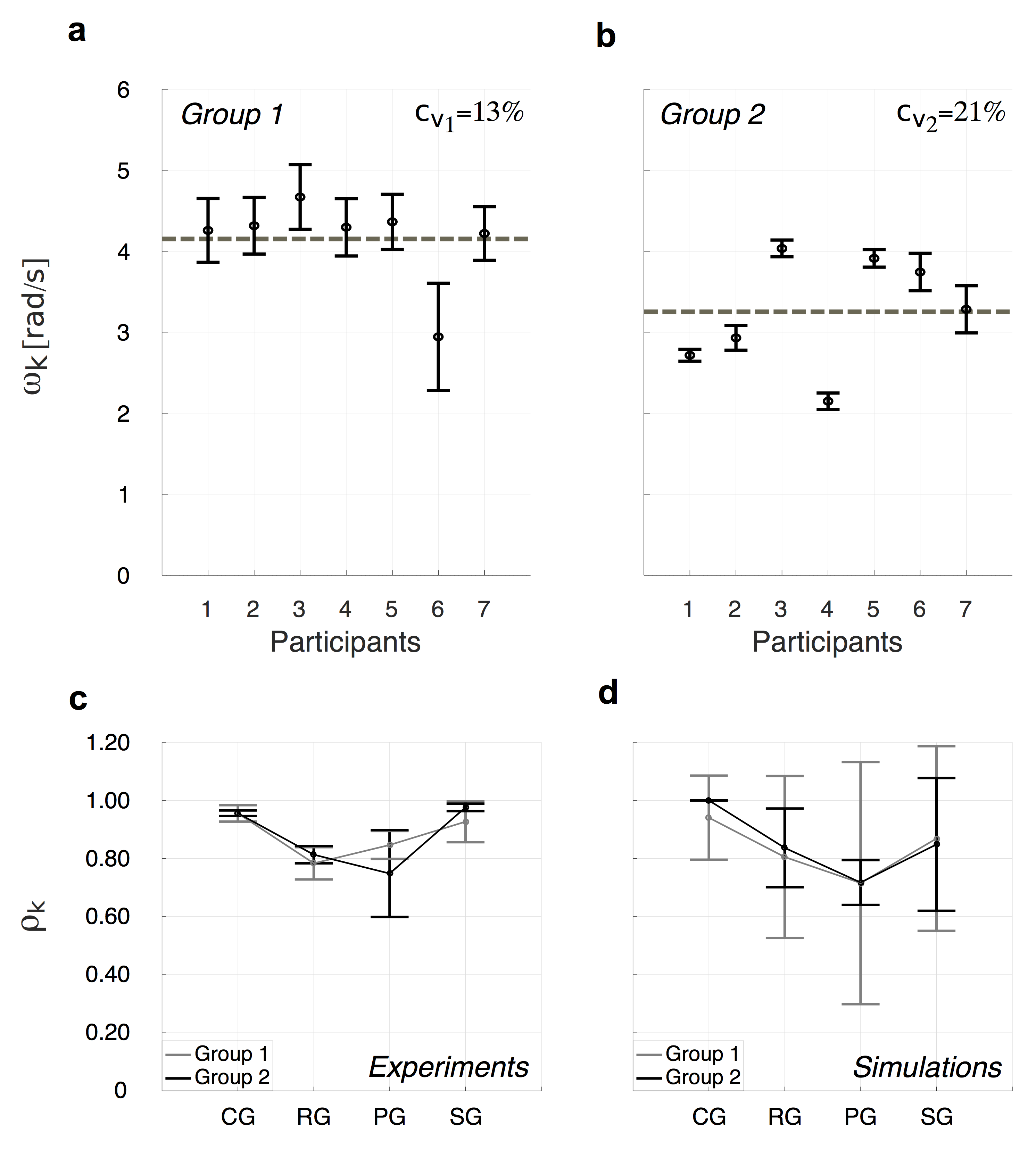}
 \caption{\textbf{Natural oscillation frequencies and individual synchronisation indices $\rho_k$ for each group and topology.} Mean (black circle) and standard deviation (black error bar) of the natural oscillation frequencies $\omega_k$ of the participants in Group 1 (\textbf{a}) and Group 2 (\textbf{b}) are presented. The frequencies of Group 2 are distributed further from their mean value averaged over the total number of players (grey dashed line) than those of Group 1, as quantified by their respective coefficient of variation $c_v$, which is equal to $c_{v_1} = 13\%$ for Group 1 and $c_{v_2} = 21\%$ for Group 2. Individual synchronisation indices are presented for experiments (\textbf{c}) and numerical simulations (\textbf{d}).  Mean values over the total number of participants are represented by circles, and standard deviations by error bars (grey for Group 1, black for Group 2). CG: Complete graph, RG: Ring graph, PG: Path graph, SG: Star graph.}
 \label{fig:Fig1}
\end{figure}

Two groups of seven players were considered, respectively named Group 1 and Group 2. Members of each group were asked to perform a simple oscillatory movement with their preferred hand and to synchronise their motion (see \emph{Methods}). The oscillations produced by each individual, when isolated from the others, had a specific natural frequency. The two groups exhibited a different level of dispersion with regards to the natural oscillation frequencies of their respective members (measured in the absence of coupling, see Section 3 of Supplementary Information), as quantified by the ensembles' coefficient of variations $c_v$ (Figs \ref{fig:Fig1}a and \ref{fig:Fig1}b). In particular, the frequencies of the players of Group 2 ($c_{v_2} = 21\%$) showed a higher dispersion than the frequencies of those of Group 1 ($c_{v_1} = 13\%$). (See also Supplementary Tables 1 and 2.)

Four different topologies of interactions were implemented through visual coupling for each group: Complete graph, Ring graph, Path graph and Star graph (Fig. \ref{fig:Fig2}). For more details on the implementation of such interaction patterns, refer to Section 1 of Supplementary Information.

\begin{figure}[!ht]
 \centering
 \includegraphics[width=1\textwidth]{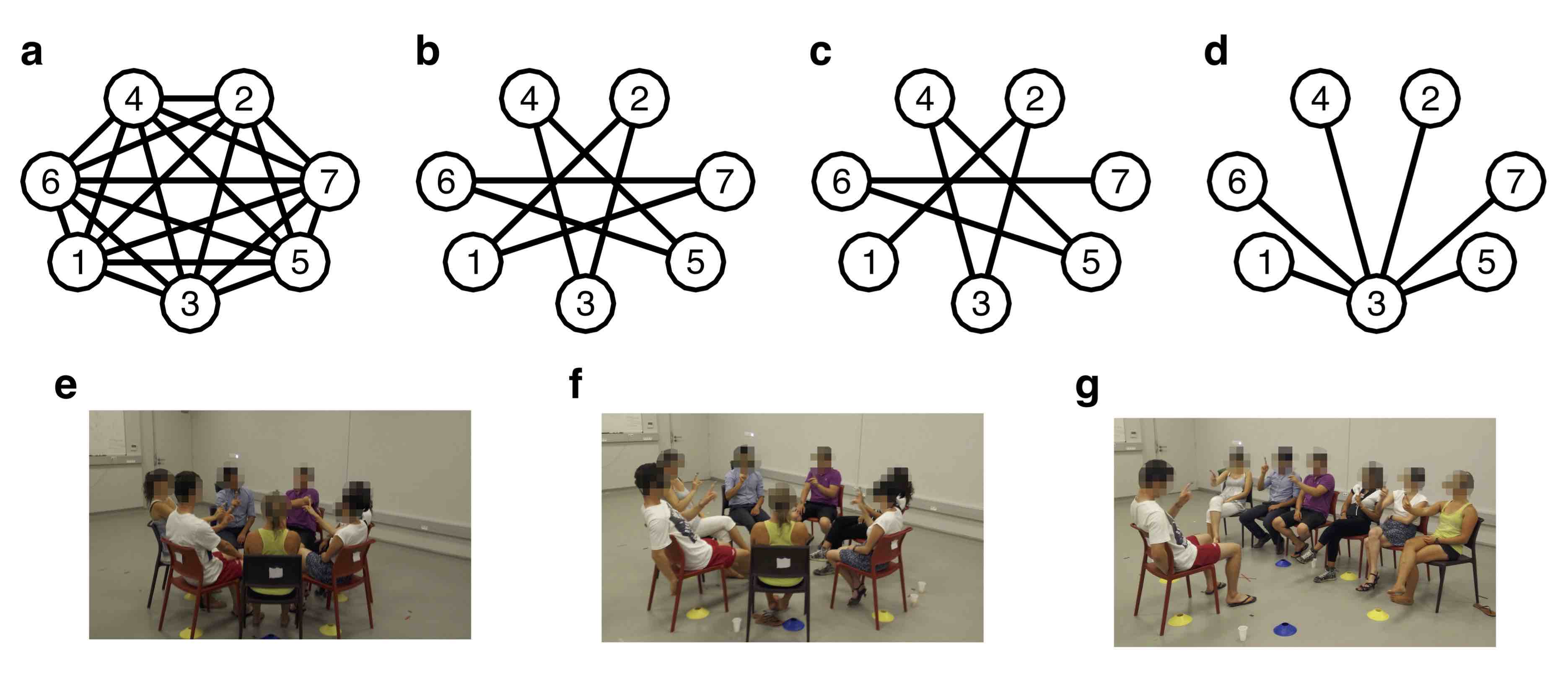}
 \caption{\textbf{Interaction patterns implemented through visual coupling in the experiments}. (\textbf{a}) Complete graph: each participant can see the movements of all the others. (\textbf{b}) Ring graph: each participant can see the movements of only her/his two \emph{partners}. (\textbf{c}) Path graph: similar to the Ring graph configuration, but agents 1 and 7, defined as \emph{external}, have only one \emph{partner} (2 and 6, respectively) and consequently are not visually coupled. (\textbf{d}) Star graph: agent 3, defined as \emph{central}, can see the movements of all the others, defined as \emph{peripheral}, who in turn see the movements of only the \emph{central} player. The other panels show the actual arrangement of the players during the experiment [(\textbf{e}) for the Complete graph, (\textbf{f}) for Ring and Path graphs, and (\textbf{g}) for the Star graph].}
 \label{fig:Fig2}
\end{figure}

A network of heterogeneous nonlinearly coupled Kuramoto oscillators \cite{K84} was employed as mathematical model to capture the relevant features observed experimentally [equation \eqref{eqn:netkoeqith}], given the oscillatory nature of the task participants were required to perform. For more details on how the parameters of such model were set, see \emph{Methods}.

\paragraph{Synchronisation levels depend on the combined action of group homogeneity and visual interactions.}

The values of the \emph{individual synchronisation indices} $\rho_k$ of the participants in the two different groups, equal to $1$ in the ideal case of player $k$ being perfectly coordinated within the ensemble and taking lower values for increasing coordination mismatches (see \emph{Methods}), were first averaged over the total number of trials for each $k$th player and for each topology, and then underwent a 2(Group) X 4(Topology) Mixed ANOVA. Their mean value and standard deviation over the total number of participants in the group are represented for each topology in Fig. \ref{fig:Fig1}c in the experimental cases, and in Fig. \ref{fig:Fig1}d for the simulations, respectively.

The ANOVA performed with the Greenhouse--Geisser correction revealed a statistically significant effect of Topology ($F(1.648,19.779)=29.447$, $p<0.01$, $\eta^2=0.710$), suggesting an advantage of both the Complete graph and the Star graph in generating higher individual synchronisation (Bonferroni post-hoc test, $p<0.01$). The Group main effect was not in itself significant ($F(1,12)=0.053$, $p=0.821$, $\eta^2=0.004$). More importantly, the significant Group X Topology interaction ($F(1.648,19.779)=3.908$, $p<0.05$, $\eta^2=0.246$) revealed that the topology effect on synchronisation was more pronounced  for the less homogenous group (Group 2), with for instance the Path graph in that group producing the lowest level of synchronisation (Bonferroni post-hoc test, $p <0.01$). For further details, see Supplementary Tables 5-7.

In short, \textbf{visual interaction between players was found to affect synchronisation indices, more so when natural individual motions differed largely from each other} \cite{STdB15}.

\begin{remark}
One could argue that the group members' plasticity, as quantified by the individual standard deviations of their natural oscillation frequencies (higher for the participants of Group 1), might be the source of differences in the overall performance. However, further numerical simulations (see Section 8.1 of Supplementary Information) confirmed that the overall frequency dispersion, rather than the intra-individual variabilities of the natural oscillation frequencies, has a significant effect on the synchronisation levels achieved by the group members. One could also argue that the difference in the natural oscillation frequencies of the participants getting disconnected in a Ring graph (to form a Path graph), or the particular member chosen as central player in a Star graph, might have a significant effect on the synchronisation levels of the ensemble. Additional numerical simulations exclude both these possibilities (see Section 8.2 and Section 8.3 of Supplementary Information for more details).
\end{remark}

\paragraph{A network of heterogeneous Kuramoto oscillators behaves like a human ensemble.}
A 2(Group) X 4(Topology) Mixed ANOVA was performed on the simulated data to evaluate the capacity of the model proposed in equation \eqref{eqn:netkoeqith} to reproduce the topology and group effects observed on the experimental human data.

The ANOVA revealed a statistically significant effect of Topology ($F(3,36)=5.946$, $p<0.01$, $\eta^2=0.331$), suggesting an advantage of the Complete graph in generating higher individual synchronisations (Bonferroni post-hoc test, $p<0.05$). The Group main effect was not significant ($F(1,12)=0.031$, $p=0.862$, $\eta^2=0.003$), and neither was the Group X Topology interaction ($F(3,36)=0.163$, $p=0.920$, $\eta^2=0.013$).
This shows that the model succeeds in replicating the statistical significant effect of Topology, with higher values of synchronisations obtained in Complete graph and Star graph, as observed experimentally. However, in its current form it fails in modulating the topology effect by variations in the group homogeneity. (See Supplementary Tables 8 and 9 for more details).

The ability of our model to capture the human synchronisation behaviour was further reinforced by the results of two Mixed ANOVAs performed with the Greenhouse--Geisser correction separately for Group 1 and for Group 2, showing no effect of Data origin (experiment vs. simulations, Group 1: $F(1,12) = 0.206$, $p=0.658$, $\eta^2=0.017$; Group 2: $F(1,12) = 0.619$, $p=0.447$, $\eta^2=0.049$), a statistical significant effect of Topology (Group 1: $F(1.523,18.272) = 5.419$, $p<0.05$, $\eta^2=0.311$; Group 2: $F(1.875,22.504) = 12.406$, $p<0.01$, $\eta^2=0.508$), and no interaction between these two factors (Group 1: $F(1.523,18.272) = 0.893$, $p=0.4$, $\eta^2=0.069$; Group 2: $F(1.875,22.504) = 1.606$, $p=0.223$, $\eta^2=0.118$). Specifically, higher synchronisations were found in the Complete graph and Star graph (Bonferroni post-hoc test, $p <0.01$). For further details, see Supplementary Tables 10-13.

Altogether, these results suggest that, \textbf{for each interaction pattern, a human ensemble and a network of Kuramoto oscillators behave similarly.} Specifically, the mathematical model proposed in equation \eqref{eqn:netkoeqith} succeeds in replicating that individual synchronisation indices $\rho_k$ depend on the particular interaction pattern implemented, as observed experimentally.

\paragraph{Effects of individual consistencies on synchronisation levels.}
Correlation analysis between individual consistencies (across trials performed in isolation) and coordination levels (obtained from group trials) ruled out the hypothesis of higher individual synchronisation indices $\rho_k$ being related to higher individual variabilities $c_v(\omega_k)$ of natural oscillation frequency (see Section 3 of Supplementary Information for more details).
Indeed, correlations between these two variables tended to be negative in the experimental data for both groups (Complete graph: $R = -0.69$, $p<0.01$; Ring graph: $R = -0.14$, $p=0.62$; Path graph: $R \simeq 0$, $p=0.83$; Star graph: $R = -0.45$, $p=0.11$), and such relations were replicated by the simulated data (Complete graph: $R = -0.88$, $p<0.01$; Ring graph: $R = -0.66$, $p<0.01$; Path graph: $R = -0.45$, $p=0.11$; Star graph: $R = -0.55$, $p<0.05$).

In short, our findings show \textbf{the existence of a negative relation between individual variabilities and synchronisation indices}, at least significantly in the Complete graph, and that the model captures such relation.

\paragraph{Visual coupling maximises synchronisation within connected dyads.}

\begin{figure}[!ht]
 \centering
 \includegraphics[width=.8\textwidth]{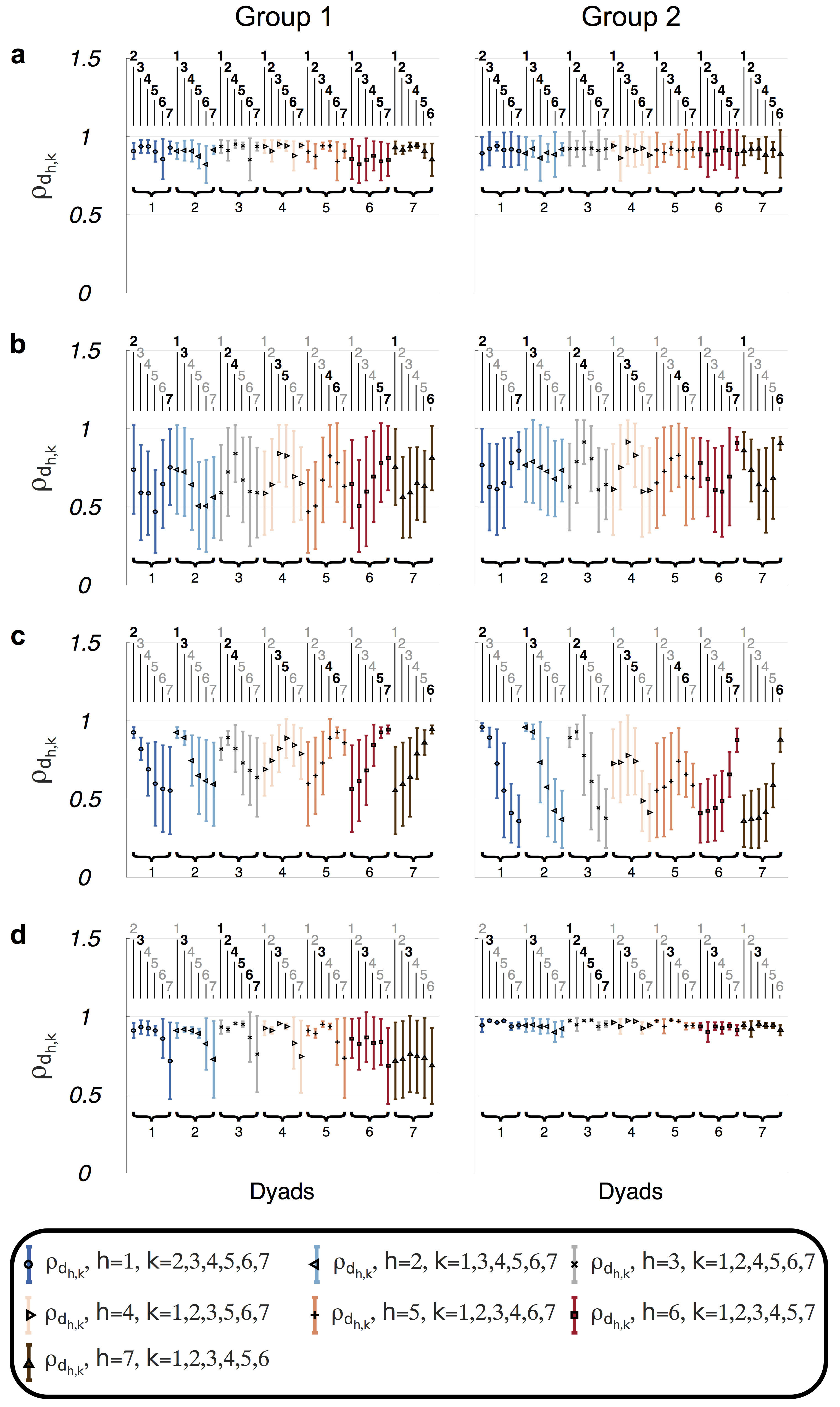}
 \caption{\textbf{Dyadic synchronisations indices observed experimentally}. Mean (symbol) and standard deviation (error bar) over the total number of trials of the \emph{dyadic synchronisation index} $\rho_{d_{h,k}}$ for players of Group 1 (left panel) and Group 2 (right panel) in Complete (\textbf{a}), Ring (\textbf{b}), Path (\textbf{c}) and Star graph (\textbf{d}) are presented (the respective interaction patterns are shown in Fig. \ref{fig:Fig2}). Different symbols and colours refer to pairs related to different players. In each panel, the black subscripts on the bottom represent $h$, whereas those on the top represent $k$ (bold black for visually connected pairs, grey for uncoupled pairs). For each participant of both groups and in all the implemented topologies, the highest mean values are obtained for the visually connected dyads.}
 \label{fig:Fig3}
\end{figure}

For each participant of both groups, in most cases ($99\%$ for Group 1 and $94\%$ for Group 2) the highest values of the \emph{dyadic synchronisation indices} $\rho_{d_{h,k}}$, defined similarly to $\rho_k$ but with respect to two generic participants $h$ and $k$ of the same group (see \emph{Methods}), were observed for the visually connected dyads, a result that was found for all topologies (Fig. \ref{fig:Fig3}).
Statistically, visually connected dyads across Group 1 and Group 2 were indeed found to exhibit higher synchronisation than non-visually coupled dyads, both in the experiments ($t(117.970) = -8.872$, $p<0.01$) and in the simulations ($t(153.326) = -6.361$, $p<0.01$). For further details, refer to Section 5 of Supplementary Information.

The trend of $\rho_{d_{h,k}}$ was thus similar in both groups, with the only exception of the Star graph in Group 1, where Player 7 failed to synchronise well with the \emph{central} node (Fig. \ref{fig:Fig3}d, left panel). In some cases, a relatively high value of \emph{group synchronisation index} $\rho_g$, equal to $1$ in the ideal case (see \emph{Methods}), coexisted with low values of $\rho_{d_{h,k}}$, as observed in the Path graph for Group 1 (Fig. \ref{fig:Fig3}c, left panel). 

This suggests that the overall group synchronisation can be high despite occasionally lower dyadic synchronisations in visually uncoupled players, as also replicated by the proposed mathematical model.
 
\paragraph{Synchronisation dynamics over time differ between group interaction patterns.}
A closer look at the dynamics of the group synchronisation index $\rho_g(t)$ over time finally revealed interesting differences between topologies. For the Complete graph, as well as to a lesser extent for the Star graph, group synchronisation was quickly reached in most trials, though with transient losses, for both Group 1 and Group 2. The scenario is different for the Ring and Path graphs, for which $\rho_g(t)$ exhibited no clear shift between transient (low time-varying) and steady state (high constant) values, but a continuous variation in its values for both human ensembles. The mathematical model proposed in equation \eqref{eqn:netkoeqith} reproduces such feature, as well as replicates the qualitative trend of $\rho_g(t)$ over time, for each group and topology (see Section 7 of Supplementary Information for further details).

These results suggest that players are not always able to maintain a high level of synchronisation over time once it is reached, particularly for the Ring and Path graphs, as also captured by the proposed mathematical model.

\clearpage

\section*{Discussion}
In this research, we studied synchronisation in multiplayer human ensembles. We asked two groups of seven volunteers to carry out a joint task in which they had to generate and synchronise an oscillatory motion of their preferred hand. 
We found that the coordination level of the participants in the two different ensembles varied over different topologies, and that such variations were more significant in the group characterised by a higher heterogeneity in terms of the natural oscillation frequency of its members.

Specifically, we observed that the synchronisation levels in Complete and Star graphs were higher than those in Path and Ring graphs, thus revealing, for the first time in a human ensemble, the key theoretical finding in \emph{Network Science} \cite{AB02,NBW11} that synchronisation depends on the structure of the interconnections between agents in a network \cite{ADKMZ08,AGGLB11,MP04,SRP06,AC16}.
In addition, we observed that the synchronisation level for a given topology was quantitatively different for the two groups characterised by a different dispersion of the agents' frequencies. In the particular case of the Path graph, we found that the group whose members had natural frequencies closer to each other (Group 1) synchronised better. This extends to multiplayer scenarios the results of \cite{ZWP16}, which showed that greater interpersonal synchrony in musical duo performances is achieved when the endogenous rhythms of two pianists are closer to each other.

Furthermore, we observed that individual consistency in the intrinsic oscillation frequency tends to enhance synchronisation, particularly in the Complete graph, and that players are not always able to maintain a high level of synchronisation over time, particularly for the Ring and Path graphs.

Note that there are other differences between the groups that might have an effect on the synchronisation levels observed experimentally (e.g., sex, weight, size, education, and so forth).
In general, social factors and personality traits may affect some of the variables defined in \emph{Methods}, and others. 
In fact, we are exploring the possibility of performing group synchronisation tasks where the same participants coordinate their motion in the presence as well as in the absence of social interaction, by means of a computer-based set-up recently proposed in \cite{ALFdB16}. Our preliminary results show that, also in the absence of social interaction, the topological structure of the interconnections among the groups members does have a statistically significant effect on their synchronisation levels.

Despite the incredible complexity of such human social interactions, we found that a rather simple mathematical model of coupled Kuramoto oscillators was able to capture most of the features observed experimentally. The availability of a mathematical description of the players' dynamics can be instrumental for designing better architectures driving virtual agents (e.g., robots, computer avatars) to coordinate their motion within groups of humans \cite{BSMGC16,IRR16,AAZdB16,ZASTdB17}, as well as for predicting the coupling strength needed to restore synchronisation based on initial knowledge of individual consistency, group variance and topology.

Even though we studied a specific laboratory-oriented joint task, our approach reveals general principles behind the emergence of movement coordination in human groups that can relate to a large variety of contexts.
A specific example where our results find confirmation is the coordination level, measured with the same \emph{group synchronisation index}, of players in a football team. As shown in \cite{RTAR13}, this index depends on the defensive playing method, giving rise to different interaction patterns among the players, as well as on the players' dynamics when considering different teams. A similar tendency was found for the synchronisation of people dancing during a club music set, which was found to depend on the features of the songs being played \cite{EBQBM16}.

More generally, our study provides a criterion to determine the best players' arrangement in multi-agent scenarios (in terms of their individual behaviour), and to designate the most appropriate interconnections among them (structure of their interactions) in order to optimise coordination when required.
This is the case for instance in music and sport, where achieving a high level of coordination is indeed a matter of crucial importance.

In music, the quality of the performance in an orchestra is related to the musicians playing in synchrony \cite{VDBCF16}. During an orchestral show (Star graph, with the central node being the conductor) the ensemble composition, in addition to classic orchestra rules, can be decided according to group heterogeneity and individual consistencies.

In collective sports, the overall performance can be improved when participants coordinate their movements \cite{RTAR13}.
For instance, in team rowing (Path graph), it is important to decide who is sitting behind whom in order to maximise group homogeneity, hence synchronisation. 
In group ice-skating, where usually athletes split into sub-groups while performing, choosing the right composition of each subset based on individual dynamics could help increase the overall coordination.
In synchronised swimming (Ring and all-to-all graphs), our findings can provide useful hints to adapt the choreographic sequence to the type of visual coupling available in these graphs.
In recreational activities (e.g., our social Sunday jogging), health benefits and social affiliation might be greater when the group members synchronise their pace \cite{MRS09}.

As we further observed that high values of $\rho_g$ can coexist with low values of $\rho_{d_{h,k}}$ for some pairs of agents, a good performance in certain group activities can be achieved by increasing specific dyadic couplings.
This is, for example, the case of people performing the Mexican wave, also known as \emph{La Ola} \cite{IHV02}. The effect of a human wave travelling across the crowd can be improved by locating side by side people who are similar in their physical characteristics as well as reaction times.

\section*{Methods}
\label{sec:method}

\paragraph{Participants.}
\label{sec:participants}
A total of 14 volunteers participated in the experiments: 5 females and 9 males (5 participants were left handed). The majority of the participants were graduate and PhD students from the EuroMov centre at the University of Montpellier in France. The experiments were held in two different sessions: seven participants took part in the first one and formed Group 1, the other seven participated in the second session and formed Group 2.

The study was carried out according to the principles expressed in the Declaration of Helsinki and was approved by the local ethical committee (EuroMov, University of Montpellier). All participants provided written informed consent for both study participation and publication of identifying information and images. Such consent was also approved by the ethical committee.

\paragraph{Task and procedure.}
\label{sec:taskproc}
Participants were asked to sit in a circle and move their preferred hand as smoothly as possible back and forth (i.e., away from and towards their bodies), along a direction required to be straight and parallel to the floor. Four different interaction patterns among players were implemented by asking each participant to focus their gaze on the motion of only a designated subset of other participants (for more details about the equipment employed and on how the different interaction structures were implemented see Section 1 of Supplementary Information).

\begin{itemize}
\item Complete graph (Figs \ref{fig:Fig2}a and \ref{fig:Fig2}e): participants were asked to keep their gaze focused on the middle of the circle in order to see the movements of all other participants.
\item Ring graph (Figs \ref{fig:Fig2}b and \ref{fig:Fig2}f): each player was asked to maintain in her/his field of view the hand motion of only two other players, called \emph{partners}.
\item Path graph (Figs \ref{fig:Fig2}c and \ref{fig:Fig2}f): similar to the Ring graph, but two participants, defined as \emph{external} participants, were asked to maintain in their field of view the hand motion of only one \emph{partner} (different for the two players).
\item Star graph (Figs \ref{fig:Fig2}d and \ref{fig:Fig2}g): all participants but one sat side-by-side facing the remaining participant. The former, defined as \emph{peripheral} players, were asked to focus their gaze on the motion of the latter, defined as \emph{central} player, who in turn was asked to maintain in her/his field of view the hand motion of all others.
\end{itemize}

Each group performed the experiments in two different conditions:
\begin{enumerate}
\item \emph{Eyes-closed condition}. Participants were asked to oscillate their preferred hand at their own comfortable tempo for $30$-second trials ($16$ trials for Group 1 and $10$ trials for Group 2) while keeping their eyes closed.
\item \emph{Eyes-open condition}. Participants were asked to synchronise the motion of each other's preferred hands during $30$-second trials. For each topology, $10$ trials lasting $30s$ each were performed.
\end{enumerate}

\paragraph{Data acquisition and analysis.}
\label{sec:data}
In order to detect the motion of the participants' hands, circular markers were attached on top of their index finger. Eight infrared cameras (Nexus MX13 Vicon System \copyright) were located around the experimental room to record the position of the markers.
For further details on how the experimental data was acquired and processed refer to Section 2 in Supplementary Information.

\paragraph{Synchronisation metrics.}
\label{sec:synchmet}

Let $x_k(t) \in \mathbb{R} \  \forall t \in [0,T]$ be the continuous time series representing the motion of each agent's preferred hand, with $k \in [1,N]$, where $N$ is the number of individuals and $T$ is the duration of the experiment. Let $x_k[t_i] \in \mathbb{R}$, with $k \in [1,N]$ and $ i \in [1,N_T]$, be the respective discrete time series of the $k$th agent, obtained after sampling $x_k(t)$ at time instants $t_i$, where $N_T$ is the number of time steps of duration $\Delta T := \frac{T}{N_T}$, that is the sampling period. Let $\theta_k(t) \in [-\pi,\pi]$ be the phase of the $k$th agent, which can be estimated by making use of the Hilbert transform of the signal $x_k(t)$ \cite{KCRPM08}. The \emph{cluster phase} or \emph{Kuramoto order parameter} is defined, both in its complex form $q'(t) \in \mathbb{C}$ and in its real form $q(t) \in [-\pi,\pi]$ as

\begin{equation}
q'(t) :=  \frac{1}{N} \sum_{k=1}^{N} e^{  j \theta_k(t)  }, \qquad q(t) := {\rm atan2} \left(\Im(q'(t)),\Re(q'(t))  \right) 
\end{equation}
which can be regarded as the average phase of the group at time $t$. 

Let $\phi_k(t) := \theta_k(t) - q(t) \in [-\pi,\pi]$ be the relative phase between the $k$th participant and the group phase at time $t$.
The relative phase between the $k$th participant and the group averaged over the time interval $[0,T]$ is defined, both in its complex form $\bar{\phi}'_k \in \mathbb{C}$ and in its real form $\bar{\phi}_k \in [-\pi,\pi]$ as

\begin{equation}
\bar{\phi}'_k := \frac{1}{T} \int_{0}^{T} e^{  j \phi_k(t)  } \ dt \simeq \frac{1}{N_T} \sum_{i=1}^{N_T} e^{  j \phi_k[t_i] },\qquad  \bar{\phi}_k := {\rm atan2} \left( \Im(\bar{\phi}'_k), \Re(\bar{\phi}'_k) \right)
\end{equation}
 
In order to quantify the degree of synchronisation of the $k$th participant with respect to the group, the following parameter

\begin{equation}
\rho_k := |\bar{\phi}'_k| \quad \in [0,1]
\end{equation}
is defined as the \emph{individual synchronisation index}: the closer $\rho_k$ is to 1, the smaller the average phase mismatch between agent $k$ and the group over the whole duration $T$ of the experiment.
 
Similarly, in order to quantify the synchronisation level of the entire group at time $t$, the following parameter

\begin{equation}
\label{eqn:r2}
\rho_g(t) := \frac{1}{N} \left | \sum_{k=1}^{N} e^{  j \left( \phi_k(t)- \bar{\phi}_k \right) } \right | \quad \in [0,1]
\end{equation}
is defined as the \emph{group synchronisation index}: the closer $\rho_g(t)$ is to 1, the smaller the average phase mismatch among the agents in the group at time $t$. Its value can be averaged over the whole time interval $[0,T]$ in order to have an estimate of the mean synchronisation level of the group during the total duration of the performance:

\begin{equation}
\label{eqn:r3}
\rho_g := \frac{1}{T} \int_{0}^{T} \rho_g(t) \ dt \simeq \frac{1}{N_T} \sum_{i=1}^{N_T} \rho_g[t_i] \quad \in [0,1]
\end{equation}

Moreover, by denoting with $\phi_{d_{h,k}}(t):=\theta_h(t)-\theta_{k}(t) \in [-\pi,\pi]$ the relative phase between two participants in the group at time $t$, it is possible to define the following parameter

\begin{equation}
\label{eqn:r4}
\rho_{d_{h,k}} := \left | \frac{1}{T} \int_{0}^{T} e^{ \{ j \phi_{d_{h,k}}(t) \} } \ dt \right | \simeq \left | \frac{1}{N_T} \sum_{i=1}^{N_T} e^{ \{ j \phi_{d_{h,k}}(t_i) \} } \right | \quad \in [0,1]
\end{equation}
as their \emph{dyadic synchronisation index}: the closer $\rho_{d_{h,k}}$ is to 1, the lower the phase mismatch between agents $h$ and $k$ over the whole trial.

\paragraph{Networks of heterogeneous Kuramoto oscillators.}
\label{sec:hetnetko}
A network of heterogeneous nonlinearly coupled Kuramoto oscillators was employed to capture the group dynamics observed experimentally \cite{K84}:

\begin{equation}
\label{eqn:netkoeqith}
\dot{\theta}_k = \omega_k + \frac{c}{N} \sum_{h=1}^{N} a_{kh} \sin \left( \theta_h - \theta_k \right), \qquad k=1,2,\ldots,N
\end{equation}
where $\theta_k$ represents the phase of the motion of the preferred hand of the $k$th human participant in the ensemble, $\omega_k$ her/his own preferred oscillation frequency when not connected to any other agent (estimated from the eyes-closed trials), and $N$ the number of participants. Each player is modelled with a different value of $\omega_k$, thus accounting for human-to-human variability, and is affected by the interaction with her/his neighbours modelled by the second term in the right hand side of equation \eqref{eqn:netkoeqith}. Speecifically, $a_{kh}=1$ if there is a connection between players $k$ and $h$ (they are looking at each other in the eyes-open trials), while $a_{kh}=0$ if there is not.

\begin{figure}[h]
 \centering
 \includegraphics[width=1\textwidth]{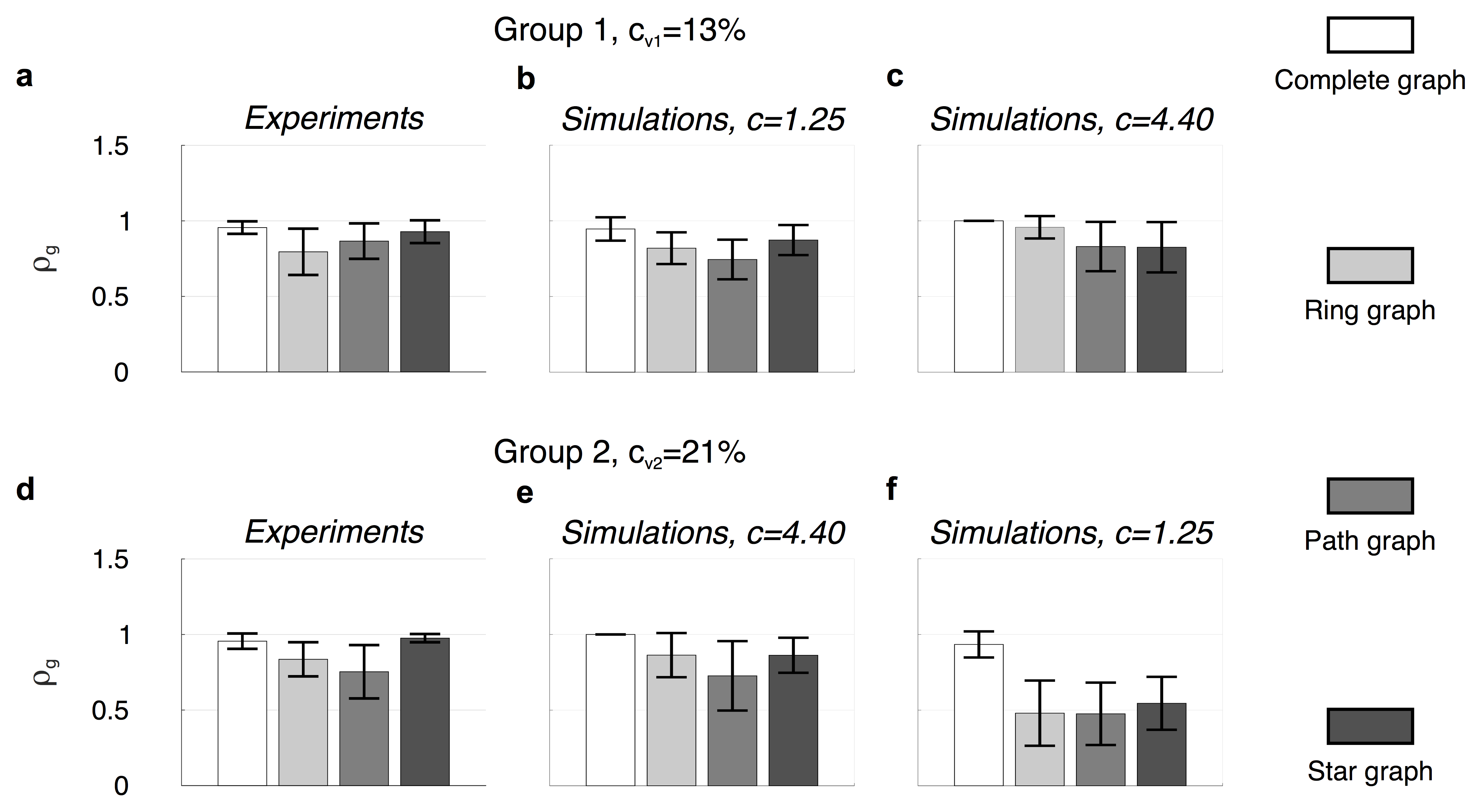}
 \caption{\textbf{Group synchronisation indices for each group and topology.} The height of each bar represents the mean value over time of $\rho_g(t)$ averaged over the total number of eyes-open trials, with different scales of grey referring to different topologies, whilst the black error bar represents its averaged standard deviation. The group synchronisation indices obtained experimentally across the four implemented topologies for Group 1 (\textbf{a}), which is characterised by a lower coefficient of variation $c_{v_1}=13\%$, are captured well numerically when $c = 1.25$ (\textbf{b}), while they are not when $c = 4.40$ (\textbf{c}). Analogously for Group 2 (\textbf{d}), which is characterised by a higher coefficient of variation $c_{v_2}=21\%$, they are captured well numerically when $c = 4.40$ (\textbf{e}), while they are not when $c = 1.25$ (\textbf{f}).}
 \label{fig:Fig4}
\end{figure}

Parameter $c$, here assumed to be constant and equal for all nodes in the network, models the interaction strength among the players, i.e., the strength of their mutual visual coupling.
Such coupling strength was set for the proposed mathematical model to match the values of group synchronisation indices $\rho_g$ observed experimentally, that is $c=1.25$ for Group 1 and $c=4.40$ for Group 2 (Fig. \ref{fig:Fig4}).

For both human ensembles, the group synchronisation indices observed experimentally (Figs \ref{fig:Fig4}a and \ref{fig:Fig4}d) is shown together with that obtained numerically by simulating the model proposed in equation \eqref{eqn:netkoeqith} with two different values of coupling strength $c$. It is possible to appreciate how a lower (higher) value of $c$ in the model reproduces well the experimental observations in the case of lower (higher) heterogeneity in the natural oscillation frequencies of the agents (as quantified by the coefficient of variation $c_v$, Figs \ref{fig:Fig4}b and \ref{fig:Fig4}e). On the other hand, experiments are not well reproduced when:
\begin{itemize}
\item the natural oscillation frequencies of the agents are close to each other and the coupling strength is too high ($c=4.40$ in Group 1, the coordination level in Complete graph and Star graph should be higher than that in Ring graph and Path graph, Fig. \ref{fig:Fig4}c);
\item the natural oscillation frequencies of the agents are far from each other and the coupling strength is too low ($c=1.25$ in Group 2, the coordination level in Ring, Path and Star graph is not comparable with that obtained experimentally, Fig. \ref{fig:Fig4}f).
\end{itemize}

As expected from theory \cite{dLdBL15}, in order to reproduce the experiments (see also Supplementary Tables 3 and 4 for further details) the coupling strength $c$ among the nodes in the model needs to take higher values for higher dispersions of the oscillation frequencies.
Further information on how the model was initialised and parameterised can be found in Section 3 of Supplementary Information.

\paragraph*{Acknowledgements.}
The authors wish to acknowledge support from the European Project AlterEgo FP7 ICT 2.9 -- Cognitive Sciences and Robotics, Grant Number 600610.
The authors wish to thank Simon Pla (Montpellier University) for the help provided in running some of the experiments, Mathieu Gueugnon (Montpellier University) for providing some of the Matlab code necessary to analyse the data, Dr. Chris Kent (the University of Bristol) for advice on how to analyse some of the data, and Prof. John Hogan (the University of Bristol), Dr. Daniel Alberto Burbano Lombana (the University of Naples Federico II) and Dr. Filippo Menolascina (the University of Edinburgh) for stimulating discussions and critical reading of the manuscript.

\paragraph*{Author contributions.}
Conceived and designed the experiments: FA, RS, BB, MdB. Performed the experiments: FA, RS. Analysed the data: FA, GF. Contributed reagents/materials/analysis tools: FA, GF, RS. Wrote the paper: FA, GF, BB, MdB.

\paragraph*{Competing financial interest.}
The authors declare no competing financial interest.

\clearpage

\section*{SUPPLEMENTARY INFORMATION}

\section{Experimental protocol}
\begin{figure}[h]
\centering
 \includegraphics[width=1\textwidth]{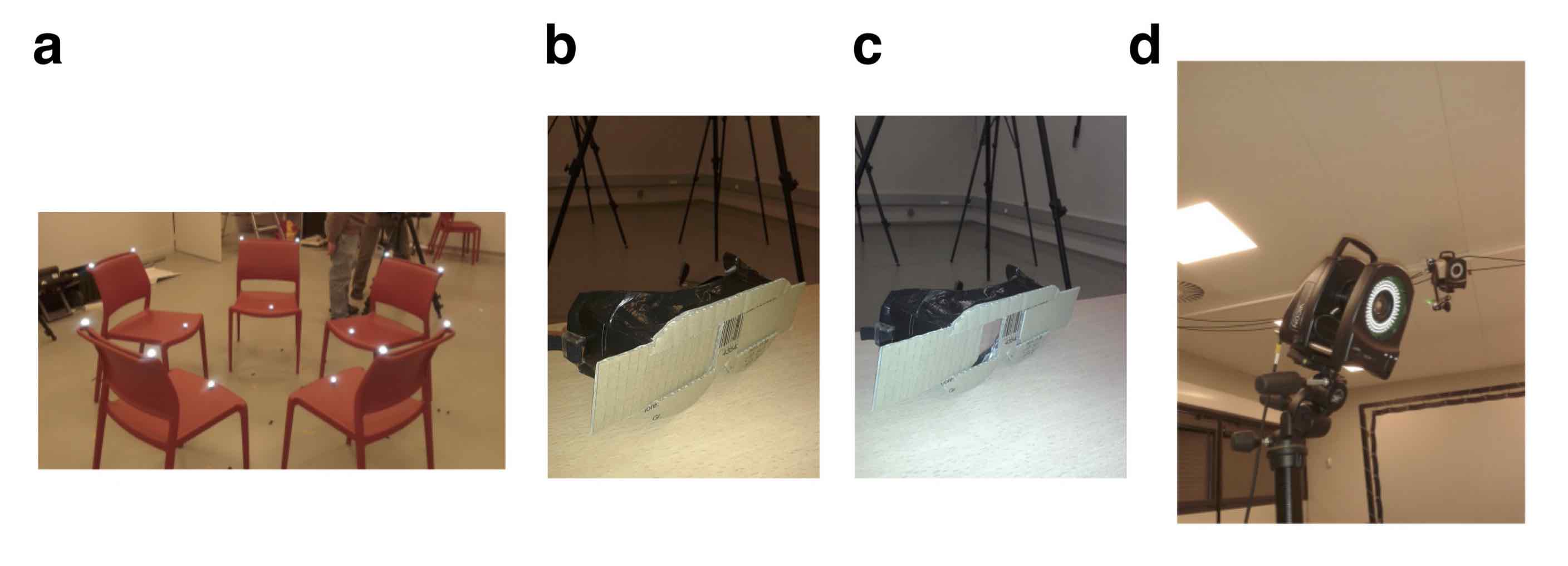}
 \caption{{\bf Experimental platform used for the experiments.} (\textbf{a}) Plastic chairs with position markers. (\textbf{b}) Shut goggles: players are deprived of their own sight. (\textbf{c}) Open goggles: by appropriately sliding the mobile cardboard on the fixed one, players can modify their own field of view. (\textbf{d}) Eight cameras are employed to track the position of the hand of each player.}
 \label{fig:S1_fig1}
\end{figure}

Participants were asked to sit in a circle on plastic chairs (Supplementary Fig. \ref{fig:S1_fig1}a) and move their preferred hand as smoothly as possible back and forth (that is, away from and back to their torso) on a plane parallel to the floor, along a direction required to be straight.

Different interaction patterns (also referred to in the text as interaction structure or topology) were implemented by asking each of the players to take into consideration the motion of only a designated subset of the other players. In order to physically implement these different interconnections, the field of view (FOV) of some players was reduced by means of \emph{ad hoc} goggles. In particular, black duct tape was wrapped around such goggles in order to mask the peripheral FOV of the players on both sides. In addition, two cardboards, one of which was mobile, were appropriately glued on the goggles so as to restrict the FOV angle of each player by adjusting the position of the sliding cardboard (Supplementary Figs \ref{fig:S1_fig1}b and \ref{fig:S1_fig1}c). In such a way it was possible to implement different visual pairings among the players.

\begin{itemize}
\item {\bf Complete graph} (Figs 2a,e in the main text): participants sat in a circle facing each other without wearing the goggles. They were asked to keep their gaze focused on the middle of the circle in order to see the movements of all the others.
\item {\bf Ring graph} (Figs 2b,f in the main text): participants sat in a circle facing each other while wearing the goggles. Each player  was asked to see the hand motion of only two others, called \emph{partners}. The goggles allowed participants to focus their gaze on the motion of their only two designated \emph{partners}.
\item {\bf Path graph} (Figs 2c,f in the main text): similar to the Ring graph configuration, but two participants, defined as \emph{external}, were asked to see the hand motion of only one \emph{partner} (not the same). This was realised by removing the visual coupling between any pair of participants in the Ring graph configuration.
\item {\bf Star graph} (Figs 2d,g in the main text): all participants but one sat side by side facing the remaining participant while wearing the goggles. The former, defined as \emph{peripheral} players, were asked to focus their gaze on the motion of the latter, defined as \emph{central} player, who in turn was asked to see the hand motion of all the others.
\end{itemize}

\section{Data analysis}

In order to detect and analyse the three-dimensional position of the participants' hands, eight infrared cameras (Nexus MX13 Vicon System \copyright) were located around the experimental room (Supplementary Fig. \ref{fig:S1_fig1}d). Despite each of them achieving full frame synchronisation up to $350Hz$, data was recorded with a sampling frequency of $100Hz$, with an estimated error of $0.01mm$ for each coordinate.

\begin{figure}[h]
\centering
 \includegraphics[width=.55\textwidth]{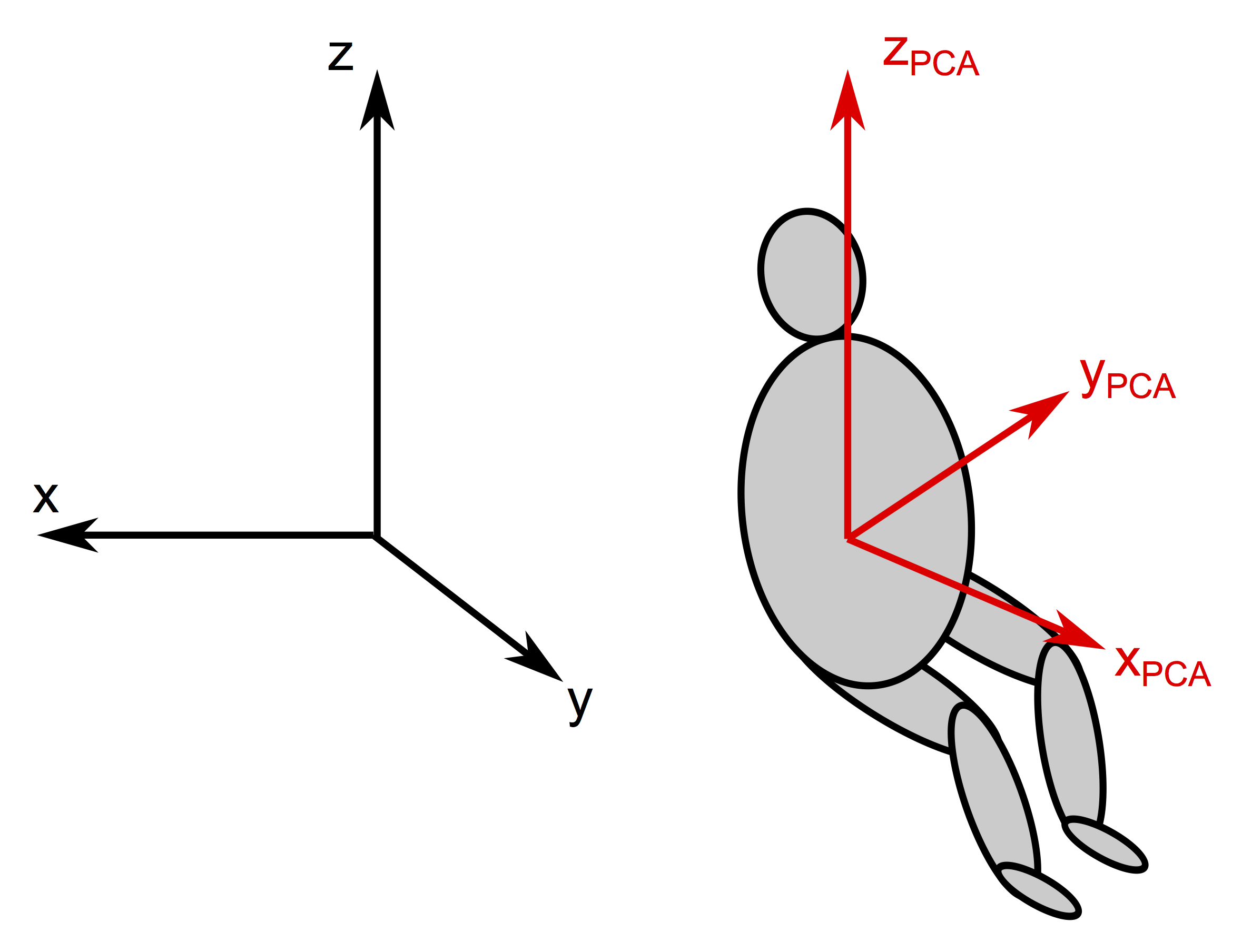}
 \caption{{\bf Cartesian frame of reference $(x,y,z)$ and principal components $(x_{PCA},y_{PCA},z_{PCA})$.} The axes $x$ and $y$ lie on a plane which is parallel to the ground, while the $z$ axis is orthogonal to it. The axes $x_{PCA}$, $y_{PCA}$ and $z_{PCA}$ individuate the principal components: $x_{PCA}$ is the direction where most of the movement takes place.}
 \label{fig:S1_fig2}
\end{figure}

In order for the cameras to detect the position of each player's hand, circular markers were attached on top of their index fingers; such positions were provided as triplets of $\left( x,y,z \right)$ coordinates in a Cartesian frame of reference (Supplementary Fig. \ref{fig:S1_fig2}, black axes). In rare occasions ($0.54 \%$ of the total number of data points for Group 1, never for Group 2) it was necessary to deprive these trajectories of possible undesired spikes caused by the cameras not being able to appropriately detect the position of the markers for the whole duration of the trial. As for Group 1, spikes were found in:
\begin{itemize}
\item 1 trajectory of Player 2 in the Ring graph topology;
\item 2 trajectories of Player 2 in the Path graph topology;
\item 4 trajectories of Player 2 and 8 trajectories of Player 7 in the Star graph topology.
\end{itemize}

After removing possible spikes, classical interpolation techniques were used to fill the gap previously occupied by the spikes themselves. Besides, since the players' positions were provided as triplets of $\left( x,y,z \right)$ coordinates but essentially the motion of each player could be described as a one-dimensional movement, it was necessary to perform \emph{principal component analysis} (PCA) on the collected trajectories to find such direction, which turns out to correspond to the $x_{PCA}$ axis (Supplementary Fig. \ref{fig:S1_fig2}, in red).

\begin{figure}[h]
\centering
 \includegraphics[width=.8\textwidth]{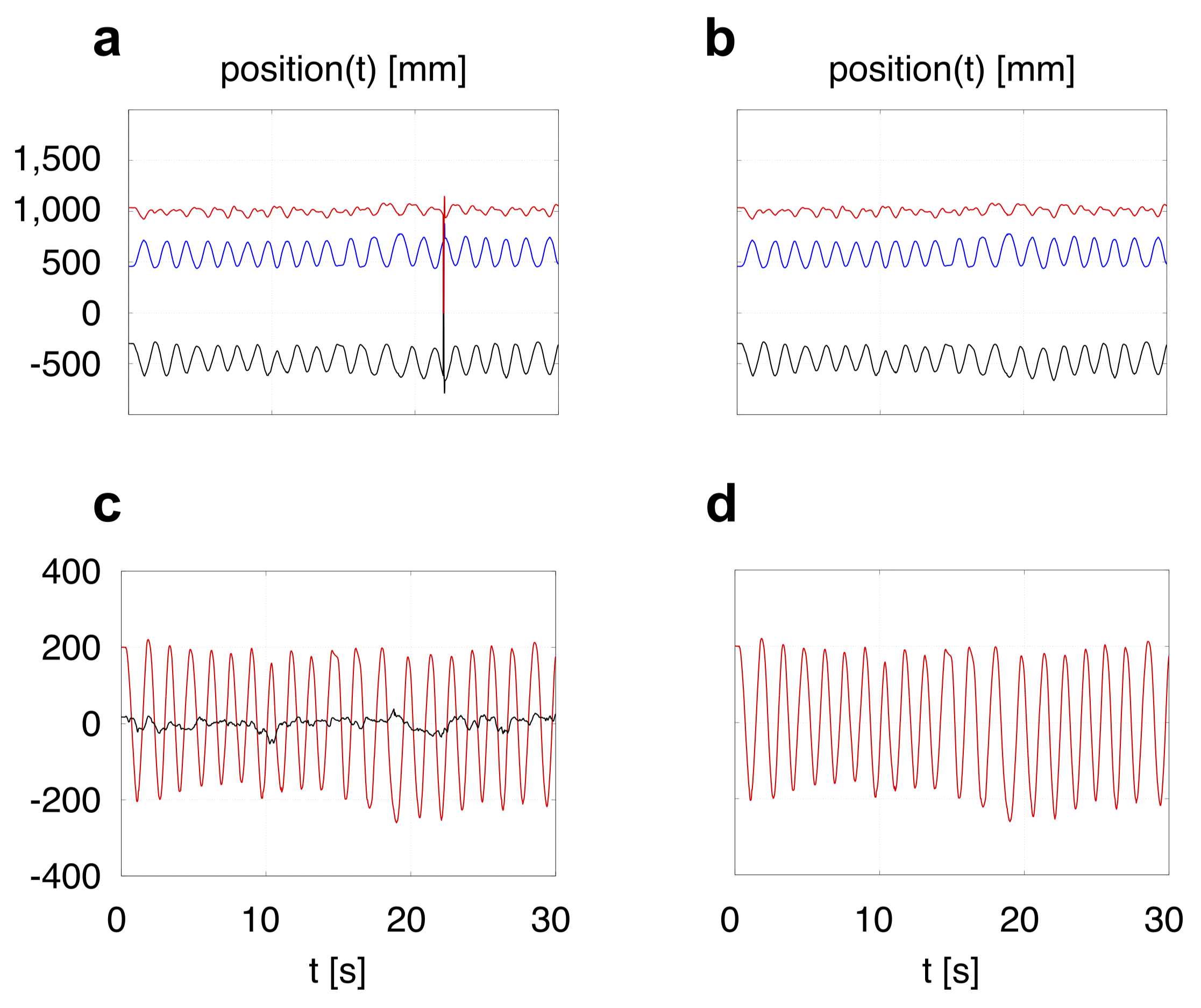}
 \caption{{\bf Pre-analysis for the hand motion of a given participant.} (\textbf{a}) Original trajectories defined by $(x,y,z)$, respectively represented in blue, black and red, after recording the hand motion of each participant. (\textbf{b}) Trajectories $(x,y,z)$, respectively represented in blue, black and red, after removing the spikes. (\textbf{c}) Trajectories $(x_{PCA},y_{PCA})$, respectively represented in red and black, after PCA analysis applied onto $x$ and $y$. (\textbf{d}) One-dimensional trajectory, defined by $x_{PCA}$ (in red), used for further analysis.}
 \label{fig:S1_fig3}
\end{figure}

PCA was applied to the collected players' trajectories defined by $x$ and $y$, since the hand motion took place mostly in the $(x,y)$ plane so that the $z$-coordinate could be neglected, obtaining the principal components $x_{PCA}$ and $y_{PCA}$. Since the motion along the component $y_{PCA}$ turns out to be negligible compared to that along $x_{PCA}$, it is possible to further assume that the motion of each player is one-dimensional (Supplementary Fig. \ref{fig:S1_fig3}). For this reason, after removing possible spikes, all the data collected in the experiments underwent PCA and then only the first principal component, namely $x_{PCA}$, was considered for further analysis.

\section{Parameterisation and initialisation of the mathematical model}

Given the oscillatory nature of the task participants were required to perform, we used a network of heterogeneous nonlinearly coupled Kuramoto oscillators as mathematical model to capture the experimental observations 

\begin{equation}
\label{eqn:netkoeqith}
\dot{\theta}_k = \omega_k + \frac{c}{N} \sum_{h=1}^{N} a_{kh} \sin \left( \theta_h - \theta_k \right), \qquad k=1,2,\ldots,N
\end{equation}

The values of the players' natural oscillation frequencies $\omega_k$ were estimated by considering the $M$ eyes-closed trials ($M = 16$ for Group 1, and $M = 10$ for Group 2). Specifically, we evaluated the fundamental harmonic of each player's motion from the Fourier transform of their position trajectory, thus obtaining $M$ values for each participant. For our simulations, we assumed that each oscillation frequency $\omega_k$ of Supplementary equation \eqref{eqn:netkoeqith} was a time-varying quantity, randomly extracted from a Gaussian distribution whose mean $\mu(\omega_k)$ and standard deviation $\sigma(\omega_k)$ are evaluated from the $M$ aforementioned values collected for each human participant (see Supplementary Fig. \ref{fig:S1_fig4} and Supplementary Tables \ref{table1} and \ref{table2} for more details).
Indeed, the Kolmogorov-Smirnov test decision for the null hypothesis that the experimental data comes from a normal distribution was performed, and such test always failed to reject the null hypothesis at the $5\%$ significance level.

\begin{figure}[h]
\centering
 \includegraphics[width=1\textwidth]{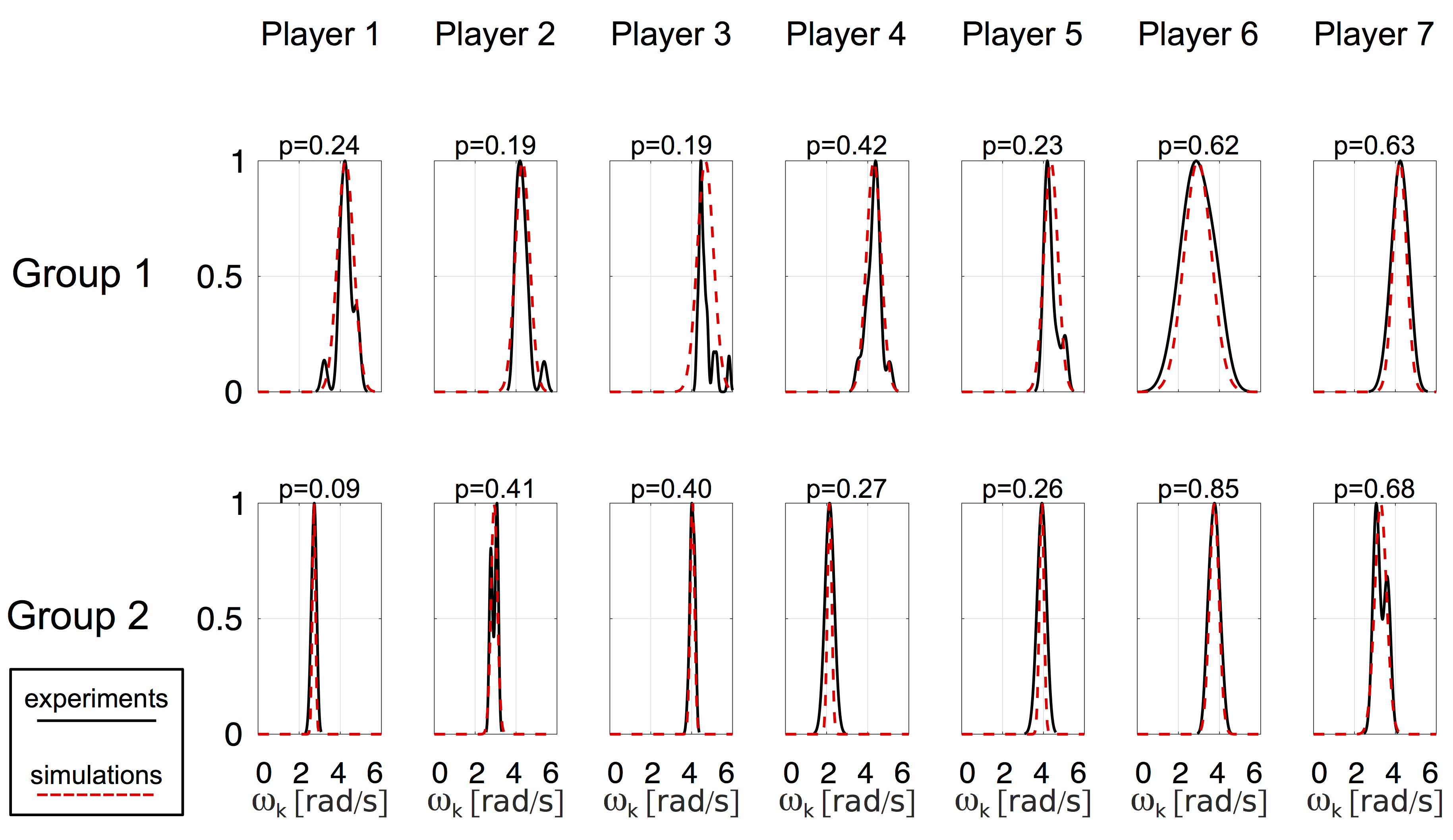}
 \caption{{\bf Probability distribution function of natural oscillation frequencies $\omega_k$.} The probability distribution functions evaluated from the $M$ values of $\omega_k$ obtained experimentally in the eyes-closed trials are represented as black solid lines, whereas the fitted normal distributions used in the numerical simulations are represented as red dashed lines.
 The null hypothesis that the experimental data comes from a normal distribution was tested. Such hypothesis could never be rejected, as specified by a $p$-value always greater than $5\%$.
 The top row refers to players of Group 1, while the bottom row to those of Group 2, whereas each column refers to a different player in the group.}
 \label{fig:S1_fig4}
\end{figure}

\begin{table}[!ht]
\caption{
{\bf Mean value and standard deviation, over the total number of eyes-closed trials, of the players' natural oscillation frequencies -- Group 1.}}
\begin{tabular}{|l|l|l|l|l|l|l|l|}
\hline
{\bf Player} & {$\mu \left(\omega_k\right)$} & {$\sigma \left(\omega_k\right)$}\\ \hline
$1$ & 4.2568 & 0.3941\\ \hline
$2$ & 4.3143 & 0.3492\\ \hline
$3$ & 4.6691 & 0.3999\\ \hline
$4$ & 4.2951 & 0.3543\\ \hline
$5$ & 4.3623 & 0.3406\\ \hline
$6$ & 2.9433 & 0.6609\\ \hline
$7$ & 4.2184 & 0.3314\\ \hline
\end{tabular}
\begin{flushleft} The mean value of the frequencies is indicated with $\mu \left(\omega_k\right)$, while their standard deviation is indicated with $\sigma \left(\omega_k\right)$, $\forall k \ \in [1,N]$. 
\end{flushleft}
\label{table1}
\end{table}

\begin{table}[!ht]
\caption{
{\bf Mean value and standard deviation, over the total number of eyes-closed trials, of the players' natural oscillation frequencies -- Group 2.}}
\begin{tabular}{|l|l|l|l|l|l|l|l|}
\hline
{\bf Player} & {$\mu \left(\omega_k\right)$} & {$\sigma \left(\omega_k\right)$}\\ \hline
$1$ & 2.7151 & 0.0741\\ \hline
$2$ & 2.9299 & 0.1525\\ \hline
$3$ & 4.0344 & 0.1035\\ \hline
$4$ & 2.1476 & 0.1023\\ \hline
$5$ & 3.9117 & 0.1085\\ \hline
$6$ & 3.7429 & 0.2309\\ \hline
$7$ & 3.2827 & 0.2911\\ \hline
\end{tabular}
\begin{flushleft} The mean value of the frequencies is indicated with $\mu \left(\omega_k\right)$, while their standard deviation is indicated with $\sigma \left(\omega_k\right)$, $\forall k \ \in [1,N]$. 
\end{flushleft}
\label{table2}
\end{table}

\begin{figure}[h]
\centering
 \includegraphics[width=1.1\textwidth]{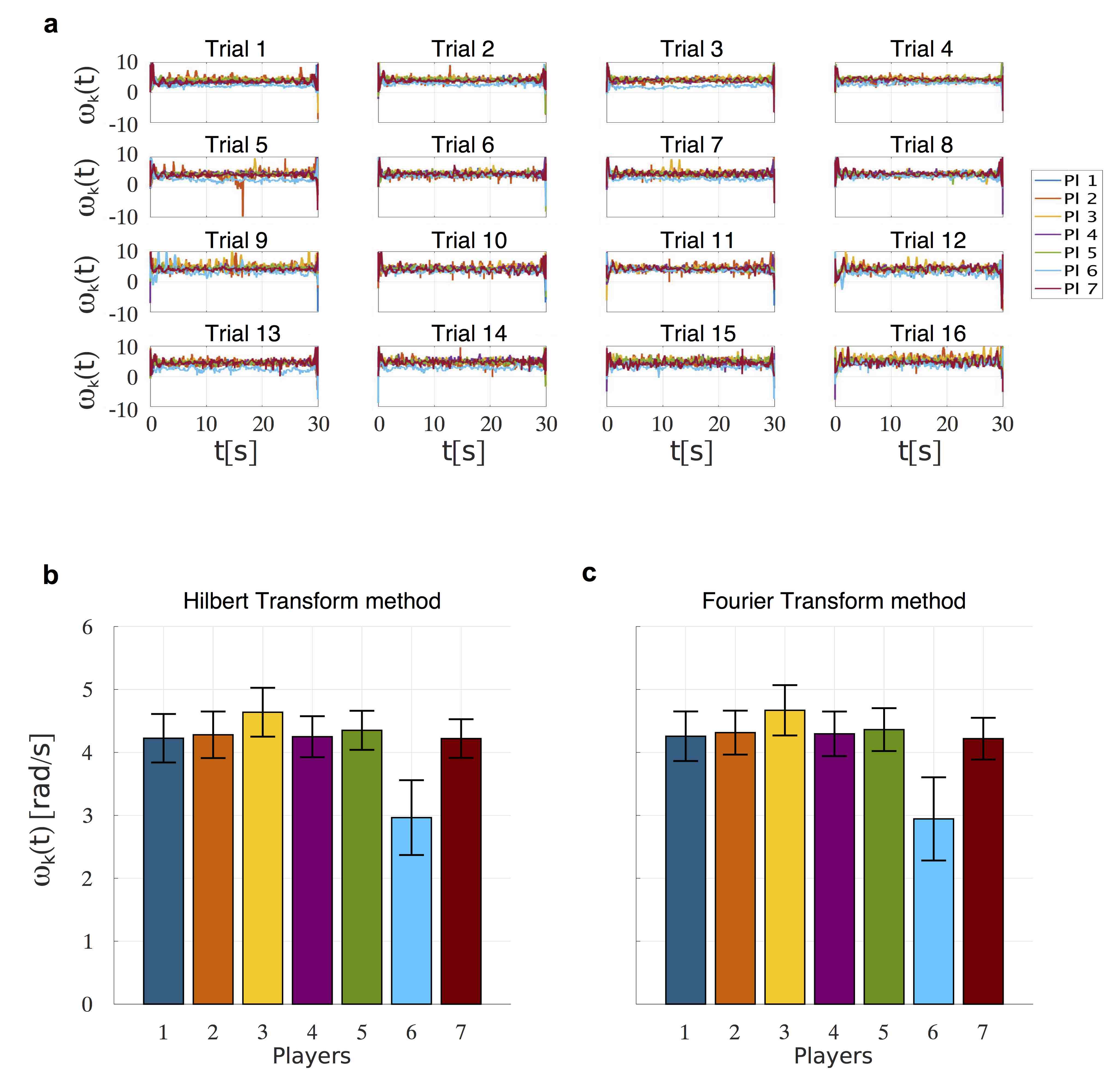}
 \caption{{\bf Natural oscillation frequencies $\omega_k$ -- Group 1.} For all the $M=16$ eyes-closed trials, the angular velocity $\omega_k(t)$ of each $k$th player, estimated through the Hilbert transform method, is a time-varying signal (\textbf{a}). Mean values (colour-coded bars) and standard deviations (black vertical bars) of the $M$ average values of $\omega_k$ over time are represented for the Hilbert transform method (\textbf{b}) and the Fourier transform method (\textbf{c}), respectively. Different colours refer to different players.}
 \label{fig:S1_fig5}
\end{figure}

\begin{figure}[h]
\centering
 \includegraphics[width=1.1\textwidth]{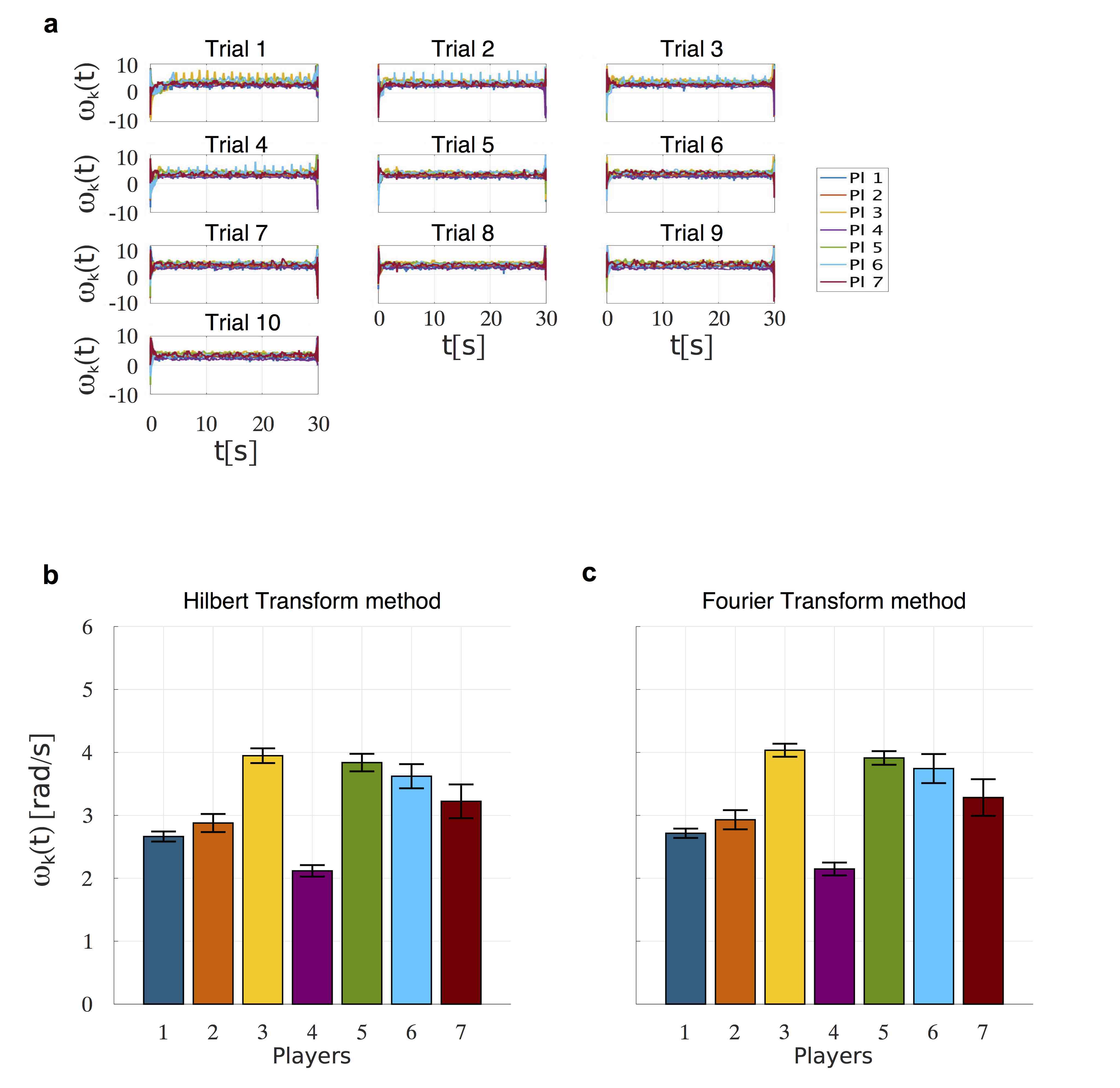}
 \caption{{\bf Natural oscillation frequencies $\omega_k$ -- Group 2.} For all the $M=10$ eyes-closed trials, the angular velocity $\omega_k(t)$ of each $k$th player, estimated through the Hilbert transform method, is a time-varying signal (\textbf{a}). Mean values (colour-coded bars) and standard deviations (black vertical bars) of the $M$ average values of $\omega_k$ over time are represented for the Hilbert transform method (\textbf{b}) and the Fourier transform method (\textbf{c}), respectively. Different colours refer to different players.}
 \label{fig:S1_fig6}
\end{figure}

Furthermore, as confirmation of the fact that players in both groups exhibited time-varying natural oscillation frequencies $\omega_k$, we also computed the Hilbert transform of each position trajectory $x_k(t)$ collected in the $M$ eyes-closed trials, then we evaluated its first time derivative, and finally observed that such derivative is a time-varying signal (Supplementary Figs \ref{fig:S1_fig5} and \ref{fig:S1_fig6}). Note that the Hilbert and the Fourier transform methods lead to consistent results for the mean values and the standard deviations of the $M$ average values of $\omega_k$ over time, respectively for each $k$th player.

By defining $\tilde{\omega} := [\mu \left(\omega_1\right) \ \mu \left(\omega_2\right) \ ... \ \mu \left(\omega_7\right)]^T \in \mathbb{R}^7$, it is possible to obtain the coefficient of variation:

\begin{equation}
c_v:=\frac{\sigma\left( \tilde{\omega} \right)}{\mu \left( \tilde{\omega} \right)}
\end{equation}
which is equal to $c_{v_1}\simeq 0.13$ for Group 1 and $c_{v_2}\simeq 0.21$ for Group 2. Such coefficient of variations quantify the overall dispersions of the natural oscillation frequencies of the players, respectively for the two groups. Analogously, it is possible to define the individual coefficient of variation:
\begin{equation}
c_v(\omega_k):=\frac{\sigma(\omega_k)}{\mu(\omega_k)}
\end{equation}
as a measure of the individual variability of the natural oscillation frequency of each $k$th player.

As for the coupling strength $c$, we found that setting the same constant value for all the topologies under investigation captures well the experimental observations (see Supplementary Section \ref{sec:gsynch} below). As for the initial values of the phases, since before starting any trial all the human players were asked to completely extend their arm so that the first movement would be pulling their arm back towards their torso from the same initial conditions, we set $\theta_k(0)=\frac{\pi}{2}$ for all the nodes, trials and topologies.

\clearpage
\section{Group synchronisation indices}
\label{sec:gsynch}

For each of the two groups we show the group synchronisation indices obtained experimentally and numerically by simulating the model proposed in Supplementary equation \eqref{eqn:netkoeqith}, with two different values of coupling strength $c$ set as described in \emph{Methods} of the main text (Supplementary Table \ref{table3} for Group 1, and Supplementary Table \ref{table4} for Group 2).

\begin{table}[!ht]
\caption{
{\bf Mean value  $\mu \left( \rho_g \right)$ and standard deviation $\sigma \left( \rho_g \right)$ over time of the \emph{group synchronisation index}, averaged over the total number of eyes-open trials -- Group 1, $c_{v_1}=13\%$.}}
\begin{tabular}{|l|l|l|l|l|l|l|l|}
\hline
{\bf Topology} & {\bf Experiments} & {\bf Simulations, $c=1.25$} & {\bf Simulations, $c=4.40$} \\ \hline
Complete graph & $0.9556 \pm 0.0414$ & $0.9462 \pm 0.0772$ & $0.9999 \pm 0.0003$ \\ \hline
Ring graph & $0.7952 \pm 0.1532$ & $0.8193 \pm 0.1048$ & $0.9575 \pm 0.0740$ \\ \hline
Path graph & $0.8661 \pm 0.1173$ & $0.7446 \pm 0.1309$ & $0.8302 \pm 0.1630$ \\ \hline
Star graph & $0.9285 \pm 0.0753$ & $0.8730 \pm 0.0993$ & $0.8255 \pm 0.1663$ \\ \hline
\end{tabular}
\begin{flushleft} This table shows $\mu \left( \rho_g \right) \pm \sigma \left( \rho_g \right)$ for both experimental and simulation results.
\end{flushleft}
\label{table3}
\end{table}

\begin{table}[!ht]
\caption{
{\bf Mean value  $\mu \left( \rho_g \right)$ and standard deviation $\sigma \left( \rho_g \right)$ over time of the \emph{group synchronisation index}, averaged over the total number of eyes-open trials -- Group 2, $c_{v_2}=21\%$.}}
\begin{tabular}{|l|l|l|l|l|l|l|l|}
\hline
{\bf Topology} & {\bf Experiments} & {\bf Simulations, $c=4.40$} & {\bf Simulations, $c=1.25$} \\ \hline
Complete graph & $0.9559 \pm 0.0508$ & $0.9999 \pm 0.0005$ & $0.9339 \pm 0.0862$ \\ \hline
Ring graph & $0.8358 \pm 0.1130$ & $0.8633 \pm 0.1460$ & $0.4799 \pm 0.2155$ \\ \hline
Path graph & $0.7534 \pm 0.1766$ & $0.7265 \pm 0.2293$ & $0.4756 \pm 0.2061$ \\ \hline
Star graph & $0.9759 \pm 0.0274$ & $0.8624 \pm 0.1158$ & $0.5450 \pm 0.1749$ \\ \hline
\end{tabular}
\begin{flushleft} This table shows $\mu \left( \rho_g \right) \pm \sigma \left( \rho_g \right)$ for both experimental and simulation results.
\end{flushleft}
\label{table4}
\end{table}


\section{Dyadic synchronisation indices}
For both Group 1 and Group 2, we show mean values and standard deviations of $\rho_{d_{h,k}}$ over the 10 eyes-open trials of each topology. In particular, if we denote with $\rho_{d_{h,k}}^{(l)}$ the value of the \emph{dyadic synchronisation index} $\rho_{d_{h,k}}$ in the $l$-th trial of a certain topology, the mean value over the total number of trials is given by

\begin{equation}
\rho_{\mu,hk} = \frac{1}{10} \sum_{l=1}^{10} \rho_{d_{h,k}}^{(l)}
\end{equation}

Similarly, the standard deviation is given by

\begin{equation}
\rho_{\sigma,hk} = \sqrt{ \frac{1}{10} \sum_{l=1}^{10} \left( \rho_{d_{h,k}}^{(l)} - \rho_{\mu,hk} \right)^2 }
\end{equation}

\begin{figure}[h]
 \centering
 \includegraphics[width=1\textwidth]{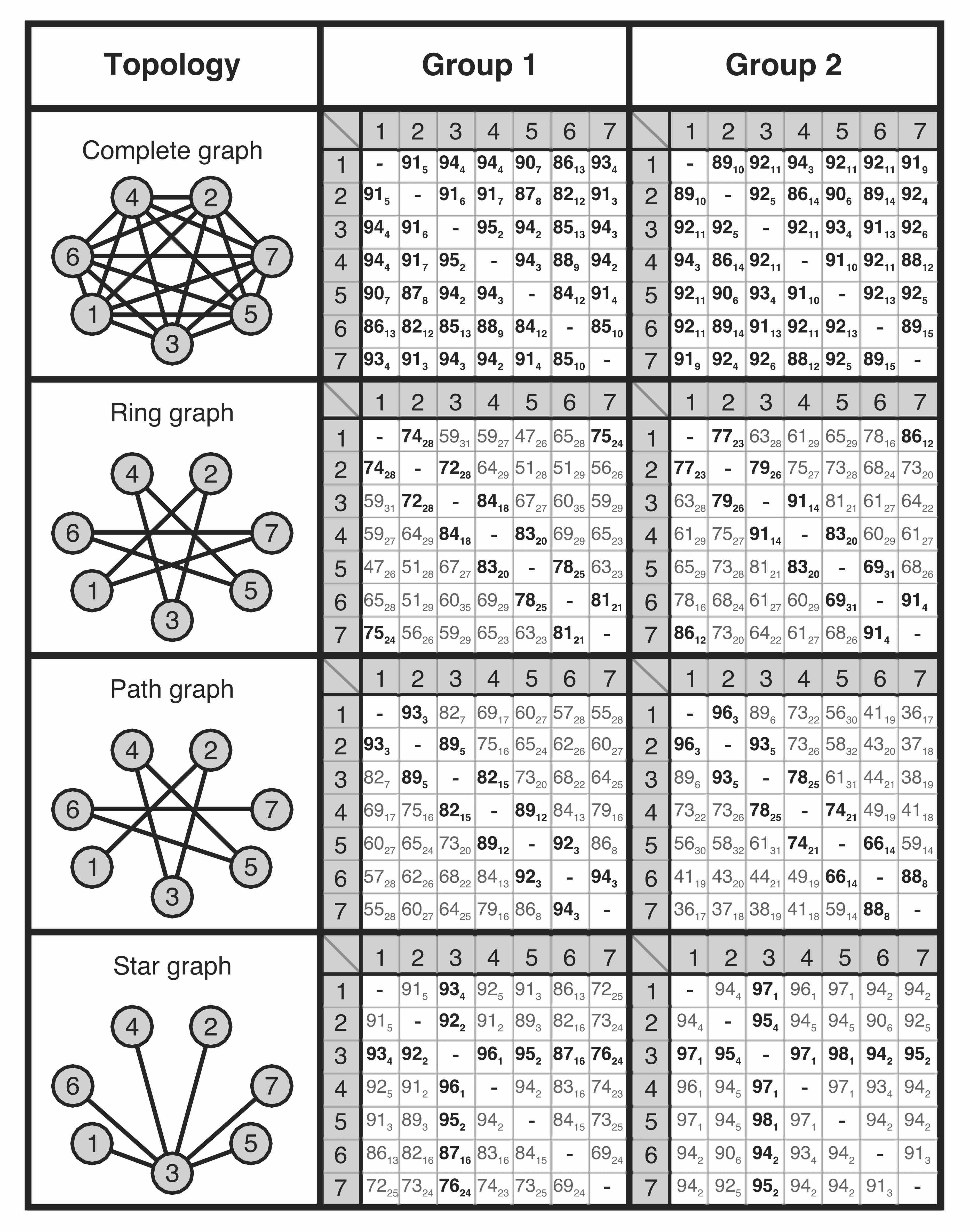}
 \caption{\textbf{Mean values and standard deviations, over the total number of eyes-open trials, of the \emph{dyadic synchronisation index} obtained experimentally.} Each row corresponds to one of the four implemented interaction patterns: topology representation (left column), indices for Group 1 (central column) and Group 2 (right column). Mean values and standard deviations (as subscripts) of $\rho_{d_{h,k}}$ are represented as percentiles for all the pairs of each topology and for both the groups, with bold values referring to pairs who were visually coupled in the experiments (i.e., there exists a link between the two agents in the respective topology representation).}
 \label{fig:S1_fig7}
\end{figure}

For both Group 1 and Group 2, mean values and standard deviations of the \emph{dyadic synchronisation index} are shown for all the pairs in the four implemented topologies (Supplementary Fig. \ref{fig:S1_fig7}). For the sake of clarity, all values of $\rho_{\mu,hk}$ and $\rho_{\sigma,hk}$ are expressed as percentiles (multiplied by $100$). In most cases ($99\%$ for Group 1 and $94\%$ for Group 2) the highest mean values of the dyadic synchronisation indices are observed for the visually connected dyads (represented in bold in Supplementary Fig. \ref{fig:S1_fig7}), meaning that players managed to maximise synchronisation with those they were visually coupled with.

In particular:

\begin{itemize}
\item Complete graph: all means $\rho_{\mu,hk}$ are higher than $0.82$ for Group 1 and $0.86$ for Group 2.
\item Ring graph: for each player of Group 1 the highest values of $\rho_{\mu,hk}$ are obtained with respect to the two agents that player was asked to be topologically connected with, whereas for each player of Group 2 at least either of the two values of $\rho_{\mu,hk}$  related to her/his \emph{partners} turns out to be the highest, and that related to the other \emph{partner} is either the second highest (nodes 2, 4 and 7), the third highest (nodes 1, 3 and 6) or the fourth highest (node 5).
\item Path graph: remarks analogous to those of the Ring graph configuration can be made. The only exception is node 4 for Group 1, where $\rho_{\mu,43}<\rho_{\mu,46}$ (however the two values are still close to each other, as $\rho_{\mu,43}=\rho_{\mu,34}=0.82$ and $\rho_{\mu,46}=\rho_{\mu,64}=0.84$), and node 3 for Group 2, where $\rho_{\mu,34}=\rho_{\mu,43}=0.78$ is lower than $\rho_{\mu,31}=\rho_{\mu,13}=0.89$. It is also worth pointing out how, for both groups, consistently with the implemented network of interactions, the mean values $\rho_{\mu,17}=\rho_{\mu,71}$ are lower than those corresponding to the Ring graph configuration, which is a consequence of removing the visual coupling between the \emph{external} agents 1 and 7.
\item Star graph: for each \emph{peripheral} player of both the groups, the highest values of $\rho_{\mu,hk}$ are obtained with respect to Player 3 (\emph{central} player).
\end{itemize}

As for the standard deviations $\rho_{\sigma,hk}$, in most cases ($86\%$ for Group 1 and $89\%$ for Group 2) the lowest values are observed for the topologically connected dyads (represented in bold in Supplementary Fig. \ref{fig:S1_fig7}), which confirms the robustness of the interactions between visually coupled pairs.

\clearpage
\section{ANOVA tables}

\begin{table}[!ht]
\centering
\caption{
{\bf 2(Group) X 4(Topology) Mixed ANOVA -- individual synchronisation indices $\rho_k$} {\bf in the experiments}}
\begin{tabular}{|l|l|l|l|l|l|l|l|}
\hline
{\bf Independent variables} & {\bf Degrees of freedom} & {\bf $F$-value} & {\bf $p$-value} & {\bf $\eta^2$} \\ \hline
Group & $(1,12)$ & $0.053$ & $0.821$ & $0.004$ \\ \hline
Topology & $(1.648,19.779)$ & $29.447$ & $\simeq 0$ & $0.710$ \\ \hline
Group * Topology & $(1.648,19.779)$ & $3.908$ & $0.044$ & $0.246$ \\ \hline
\end{tabular}
\label{table5}
\end{table}


\begin{table}[!ht]
\centering
\caption{
{\bf Post-hoc pairwise comparisons -- individual synchronisation indices $\rho_k$ in the experiments -- Group 1}}
\begin{tabular}{|l|l|l|l|l|l|l|l|}
\hline
{\bf Topologies} & Complete graph & Ring graph & Path graph & Star graph \\ \hline
Complete graph & $-$ & $\simeq 0$ & $0.146$ & $0.929$ \\ \hline
Ring graph & $\simeq 0$ & $-$ & $0.852$ & $\simeq 0$ \\ \hline
Path graph & $0.146$ & $0.852$ & $-$ & $0.564$ \\ \hline
Star graph & $0.929$ & $\simeq 0$ & $0.564$ & $-$ \\ \hline
\end{tabular}
\label{table6}
\end{table}

\begin{table}[!ht]
\centering
\caption{
{\bf Post-hoc pairwise comparisons -- individual synchronisation indices $\rho_k$ in the experiments -- Group 2}}
\begin{tabular}{|l|l|l|l|l|l|l|l|}
\hline
{\bf Topologies} & Complete graph & Ring graph & Path graph & Star graph \\ \hline
Complete graph & $-$ & $\simeq 0$ & $0.002$ & $1$ \\ \hline
Ring graph & $\simeq 0$ & $-$ & $0.792$ & $\simeq 0$ \\ \hline
Path graph & $0.002$ & $0.792$ & $-$ & $0.001$ \\ \hline
Star graph & $1$ & $\simeq 0$ & $0.001$ & $-$ \\ \hline
\end{tabular}
\label{table7}
\end{table}

\begin{table}[!ht]
\centering
\caption{
{\bf 2(Group) X 4(Topology) Mixed ANOVA -- individual synchronisation indices $\rho_k$} {\bf in the numerical simulations}}
\begin{tabular}{|l|l|l|l|l|l|l|l|}
\hline
{\bf Independent variables} & {\bf Degrees of freedom} & {\bf $F$-value} & {\bf $p$-value} & {\bf $\eta^2$} \\ \hline
Group & $(1,12)$ & $0.031$ & $0.862$ & $0.003$ \\ \hline
Topology & $(3,36)$ & $5.946$ & $\simeq 0$ & $0.331$ \\ \hline
Group * Topology & $(3,36)$ & $0.163$ & $0.920$ & $0.013$ \\ \hline
\end{tabular}
\label{table8}
\end{table}

\begin{table}[h]
\caption{
{\bf Post-hoc pairwise comparisons -- individual synchronisation indices $\rho_k$} {\bf in the numerical simulations}}
\centering 
\begin{tabular}{l l c c} 
\hline\hline 
 {\bf Topology (A)} & {\bf Topology (B)} & {\bf Mean Difference (A-B)} & {\bf $p$-value}
\\ [0.5ex]
\hline 
Complete graph & Ring graph &0.149 & 0.010 \\[-.2ex]
& Path graph &0.254 &0.017 \\[-.2ex]
& Star graph &0.112 &0.370 \\[-.2ex]
\hline 
Ring graph & Complete graph &-0.149 & 0.010 \\[-.2ex]
& Path graph &0.104 &0.795 \\[-.2ex]
& Star graph &-0.038 & 1 \\[-.2ex]
\hline 
Path graph & Complete graph &-0.254 & 0.017 \\[-.2ex]
& Ring graph &-0.104 &0.795 \\[-.2ex]
& Star graph &-0.142 &0.706 \\[-.2ex]
\hline 
Star graph & Complete graph &-0.112 & 0.370 \\[-.2ex]
& Ring graph &0.038 & 1 \\[-.2ex]
& Path graph &0.142 &0.706 \\[-.2ex]
\hline 
\end{tabular}
\label{table9}
\end{table}

\begin{table}[!ht]
\centering
\caption{
{\bf 2(Data origin) X 4(Topology) Mixed ANOVA -- individual synchronisation indices $\rho_k$ in experiments and numerical simulations -- Group 1}}
\begin{tabular}{|l|l|l|l|l|l|l|l|}
\hline
{\bf Independent variables} & {\bf Degrees of freedom} & {\bf $F$-value} & {\bf $p$-value} & {\bf $\eta^2$} \\ \hline
Data origin & $(1,12)$ & $0.206$ & $0.658$ & $0.017$ \\ \hline
Topology & $(1.523,18.272)$ & $5.419$ & $0.020$ & $0.311$ \\ \hline
Data origin * Topology & $(1.523,18.272)$ & $0.893$ & $0.400$ & $0.069$ \\ \hline
\end{tabular}
\label{table10}
\end{table}

\begin{table}[!ht]
\centering
\caption{
{\bf Main effects of Topology on individual synchronisation indices $\rho_k$ in experiments and numerical simulations -- Group 1}}
\begin{tabular}{|l|l|l|l|l|l|l|l|}
\hline
{\bf Topologies} & Complete graph & Ring graph & Path graph & Star graph \\ \hline
Complete graph & $-$ & $0.001$ & $0.174$ & $1.000$ \\ \hline
Ring graph & $0.001$ & $-$ & $1.000$ & $0.003$ \\ \hline
Path graph & $0.174$ & $1.000$ & $-$ & $0.636$ \\ \hline
Star graph & $1.000$ & $0.003$ & $0.636$ & $-$ \\ \hline
\end{tabular}
\label{table11}
\end{table}

\begin{table}[!ht]
\centering
\caption{
{\bf 2(Data origin) X 4(Topology) Mixed ANOVA -- individual synchronisation indices $\rho_k$ in experiments and numerical simulations -- Group 2}}
\begin{tabular}{|l|l|l|l|l|l|l|l|}
\hline
{\bf Independent variables} & {\bf Degrees of freedom} & {\bf $F$-value} & {\bf $p$-value} & {\bf $\eta^2$} \\ \hline
Data origin & $(1,12)$ & $0.619$ & $0.447$ & $0.049$ \\ \hline
Topology & $(1.875,22.504)$ & $12.406$ & $\simeq 0$ & $0.508$ \\ \hline
Data origin * Topology & $(1.875,22.504)$ & $1.606$ & $0.223$ & $0.118$ \\ \hline
\end{tabular}
\label{table12}
\end{table}

\begin{table}[!ht]
\centering
\caption{
{\bf Main effects of Topology on individual synchronisation indices $\rho_k$ in experiments and numerical simulations -- Group 2}}
\begin{tabular}{|l|l|l|l|l|l|l|l|}
\hline
{\bf Topologies} & Complete graph & Ring graph & Path graph & Star graph \\ \hline
Complete graph & $-$ & $\simeq 0$ & $\simeq 0$ & $0.926$ \\ \hline
Ring graph & $\simeq 0$ & $-$ & $0.390$ & $0.352$ \\ \hline
Path graph & $\simeq 0$ & $0.390$ & $-$ & $0.069$ \\ \hline
Star graph & $0.926$ & $0.352$ & $0.069$ & $-$ \\ \hline
\end{tabular}
\label{table13}
\end{table}

%

\clearpage

\section{Group synchronisation trend over time}
\begin{figure}[h]
 \centering
 \includegraphics[width=.8\textwidth]{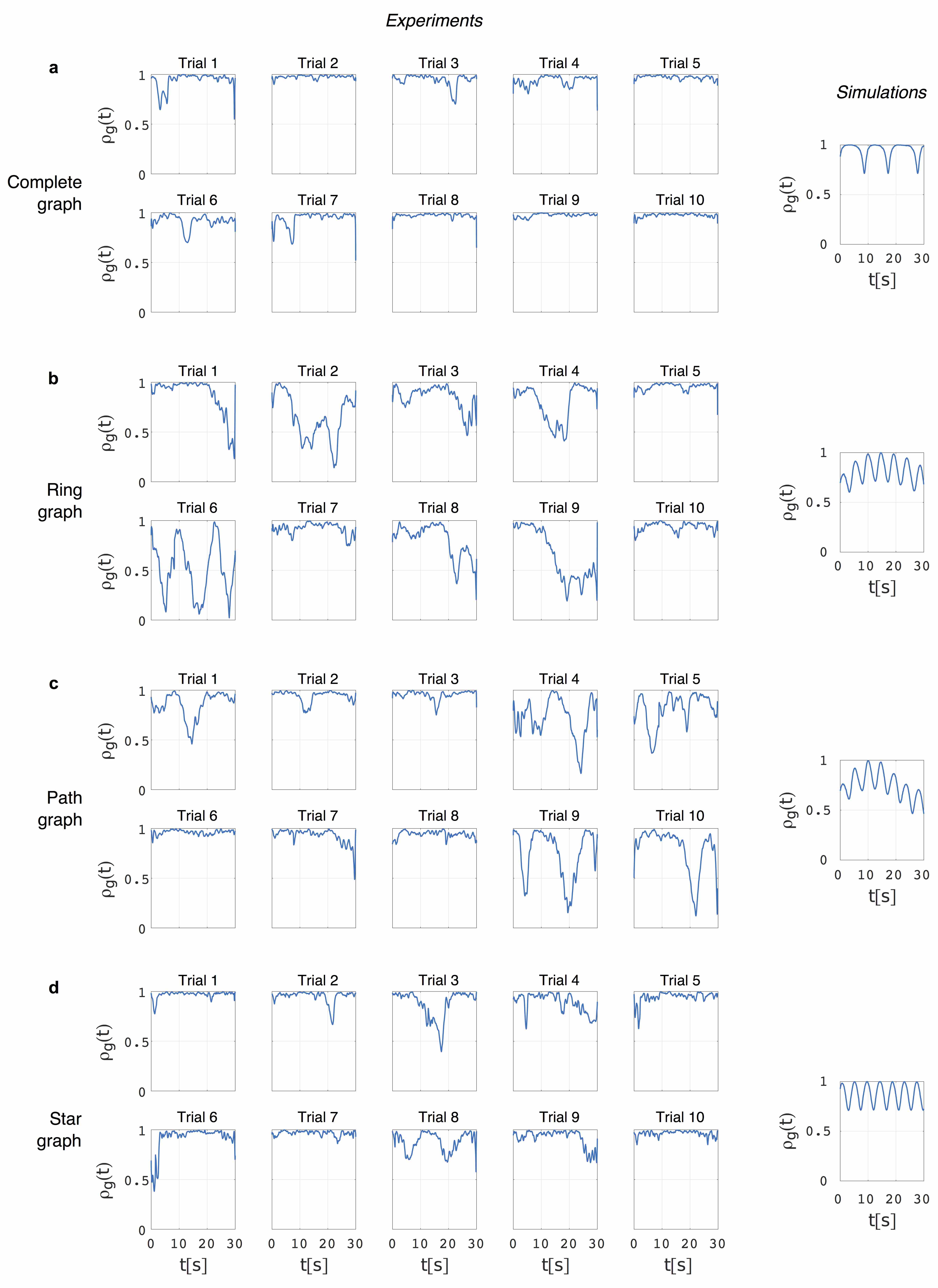}
 \caption{\textbf{Trend over time of group synchronisation index $\rho_g(t)$ for each trial and topology, both for experiments and numerical simulations -- Group 1.}  (\textbf{a}) Complete graph, (\textbf{b}) Ring graph, (\textbf{c}) Path graph, (\textbf{d}) Start Graph. For each topology, the ten panels on the left show the trend of $\rho_g(t)$ observed experimentally in the ten eyes-open trials, respectively, whereas the panel on the right shows a typical trend of $\rho_g(t)$ obtained numerically for that topology.}
 \label{fig:S1_fig8}
\end{figure}

\begin{figure}[h]
 \centering
 \includegraphics[width=.8\textwidth]{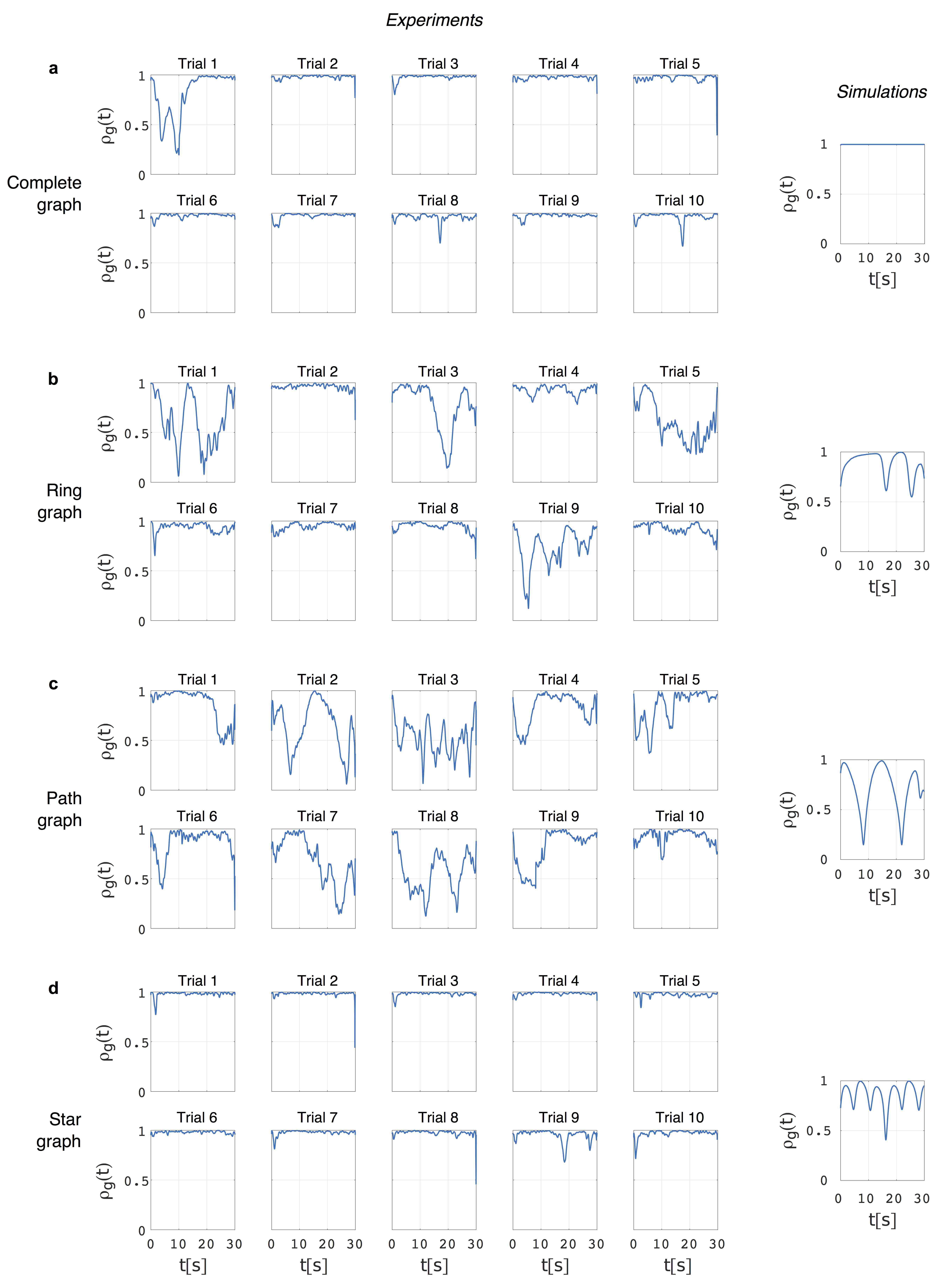}
 \caption{\textbf{Trend over time of group synchronisation index $\rho_g(t)$ for each trial and topology, both for experiments and numerical simulations -- Group 2.}  (\textbf{a}) Complete graph, (\textbf{b}) Ring graph, (\textbf{c}) Path graph, (\textbf{d}) Start Graph. For each topology, the ten panels on the left show the trend of $\rho_g(t)$ observed experimentally in the ten eyes-open trials, respectively, whereas the panel on the right shows a typical trend of $\rho_g(t)$ obtained numerically for that topology.}
 \label{fig:S1_fig9}
\end{figure}

Supplementary Fig. \ref{fig:S1_fig8} shows the trend over time of the group synchronisation index $\rho_g(t)$ observed in the experiments, as well as in the numerical simulations, for all the eyes-open trials and topologies of Group 1, whereas Supplementary Fig. \ref{fig:S1_fig9} shows that of Group 2.

It is possible to appreciate that the proposed mathematical model succeeds in replicating the feature observed experimentally that there is no clear shift between transient (low time-varying) and steady state (high constant) values for the group synchronisation index $\rho_g(t)$, and that more noticeable oscillations are observed in the Ring and Path graphs. Furthermore, note how both in the experiments and in the simulations:

\begin{itemize}
\item for Group 1 (Supplementary Fig. \ref{fig:S1_fig8}), $\rho_g(t)$ never achieves a constant value at steady state but exhibits persistent oscillations (less noticeable in the Complete graph);
\item for Group 2 (Supplementary Fig. \ref{fig:S1_fig9}), $\rho_g(t)$ exhibits higher oscillations in Ring and Path graphs, whereas lower peaks are observed in the Star graph and almost constant steady state values in the Complete graph.
\end{itemize}

\clearpage

\section{Additional model predictions}
In order to address the issues raised in \emph{Remark 1} of the main text, in this section we show further results obtained numerically by simulating a network of heterogeneous Kuramoto oscillators. Specifically, the proposed model predicts that:
\begin{enumerate}
\item unlike the intra-individual variability of oscillation frequencies $\sigma(\omega_k)$, the overall dispersion $c_{v}$ has a significant effect on the coordination levels of the group members;
\item the location of the link getting removed from a Ring graph to form a Path graph does not have a significant effect on the coordination levels of its members;
\item different selections of the central node in a Star graph do not have a significant effect on the coordination levels of its members.
\end{enumerate}

\subsection{Effects of overall dispersion and intra-individual variability of the natural oscillation frequencies}
Firstly, we considered four different groups of $N=7$ heterogeneous Kuramoto oscillators. In order to isolate the effects of the overall dispersion $c_{v}$, the individual standard deviations of the natural oscillation frequencies were not varied across nodes and groups ($\sigma(\omega_k)=0.25$ $\forall k \in [1,N]$), and all the agents in each of the four groups were connected over a Complete graph topology (so that every node is connected to all the others, and no effects of particular topological structures or symmetries are expected). Therefore, the groups differed only for the value of their overall dispersion, respectively equal to $c_{v_1}=8\%$, $c_{v_2}=17\%$, $c_{v_3}=41\%$ and $c_{v_4}=62\%$. For each group, 10 trials of duration $T=30s$ were run, with $c=1$ and initial conditions equal to $\theta_k(0)=\frac{\pi}{2}$ $\forall k \in [1,N]$.

\begin{figure}[!ht]
 \centering
 \includegraphics[width=.75\textwidth]{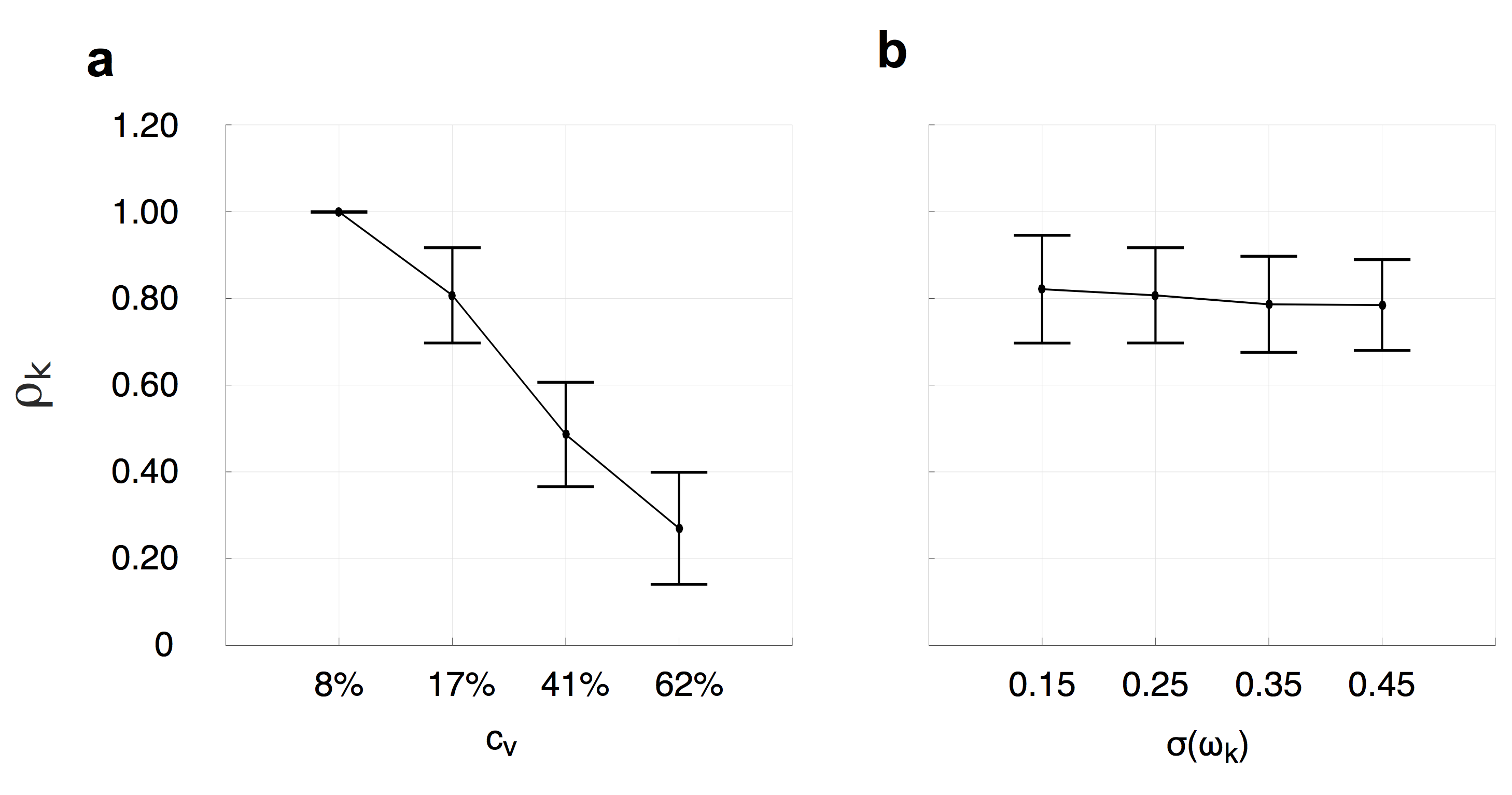}
 \caption{\textbf{Individual synchronisation indices $\rho_k$ as a function of overall dispersion $c_{v}$ and intra-individual variability $\sigma(\omega_k)$ of the natural oscillation frequencies.} Individual synchronisation indices are shown for the four groups with equal $\sigma(\omega_k)=0.25$ and different $c_{v}$ (\textbf{a}), as well as for those with equal $c_v=17\%$ and different $\sigma(\omega_k)$ (\textbf{b}). Mean values over the total number of nodes are represented by circles, and standard deviations by error bars.}
 \label{fig:S1_fig10}
\end{figure}

Supplementary Fig. \ref{fig:S1_fig10}a shows the values of the individual synchronisation indices $\rho_k$ as a function of the overall dispersion $c_{v}$. A One-way ANOVA using the Welch's test revealed a statistically significant effect of $c_{v}$ ($F(3,10)=109.345$, $p<0.001$, $\eta^2=0.896$). A post-hoc multiple comparison using the Games-Howell test is detailed in Supplementary Table \ref{table14}, showing that the differences in $\rho_k$ are statistically significant between all the groups ($p<0.05$).

\begin{table}[!ht]
\centering
\caption{
{\bf Post-hoc multiple comparisons -- individual synchronisation indices $\rho_k$ -- effects of $c_{v}$}}
\begin{tabular}{|l|l|l|l|l|l|l|l|}
\hline
$c_{v}$ & $\textbf{8\%}$ & $\textbf{17\%}$ & $\textbf{41\%}$ & $\textbf{62\%}$ \\ \hline
$\textbf{8\%}$ & $-$ & $0.014$ & $\simeq 0$ & $\simeq 0$ \\ \hline
$\textbf{17\%}$ & $0.014$ & $-$ & $0.001$ & $\simeq 0$ \\ \hline
$\textbf{41\%}$ & $\simeq 0$ & $0.001$ & $-$ & $0.031$ \\ \hline
$\textbf{62\%}$ & $\simeq 0$ & $\simeq 0$ & $0.031$ & $-$ \\ \hline
\end{tabular}
\label{table14}
\end{table}

Secondly, we considered four other different groups of $N=7$ heterogeneous Kuramoto oscillators. In order to isolate the effects of the intra-individual variability of the natural oscillation frequencies, the overall frequency dispersion was not varied over the groups ($c_v=17\%$), whereas $\sigma(\omega_k)$ was varied across them ($\sigma(\omega_k)=0.15$ for each $k$th node of the first group, $\sigma(\omega_k)=0.25$ for the second group, $\sigma(\omega_k)=0.35$ for the third group, $\sigma(\omega_k)=0.45$ for the fourth group). The other parameters were set as in previous case.

Supplementary Fig. \ref{fig:S1_fig10}b shows the values of the individual synchronisation indices $\rho_k$ as a function of the intra-individual variability $\sigma(\omega_k)$ of the natural oscillation frequencies. A One-way ANOVA revealed no statistically significant effect of $\sigma(\omega_k)$ ($F(3,24)=0.170$, $p=0.916$, $\eta^2=0.019$).

Overall, these results suggest that, unlike the intra-individual variability of oscillation frequencies $\sigma(\omega_k)$, the overall dispersion $c_{v}$ has a significant effect on the coordination levels of the group members.

\subsection{Location of the link getting removed from a Ring graph}

We then considered a network of $N=7$ heterogeneous Kuramoto oscillators connected over a Path graph topology (Fig. 2c in the main text). Seven scenarios were considered, where each scenario differs from the others in the Ring graph connection (Fig. 2b in the main text) getting removed to form the Path graph itself. Specifically, in the first scenario the connection between nodes $1$ and $2$ was removed, in the second scenario that between nodes $2$ and $3$, up to the seventh scenario where the connection between nodes $7$ and $1$ was removed. In order to isolate the effects of such choice on the coordination levels $\rho_k$, the group dispersion was not varied over the different scenarios ($c_v=17\%$), and neither were the individual variabilities across all the nodes ($\sigma(\omega_k)=0.25 \ \forall k \in [1,N]$).  The other parameters were set as in the previous cases.

A One-way ANOVA revealed that the location of the link getting removed from a Ring graph to form a Path graph, and hence the difference between the frequencies of the nodes getting disconnected, does not have a significant effect on the coordination levels of its members ($F(6,42)=0.535$, $p=0.778$, $\eta^2=0.071$). Mean values and standard deviations of the individual synchronisation index $\rho_k$ obtained in the seven different scenarios here considered are shown in Supplementary Fig. \ref{fig:S1_fig11}.

\begin{figure}[!ht]
 \centering
 \includegraphics[width=.75\textwidth]{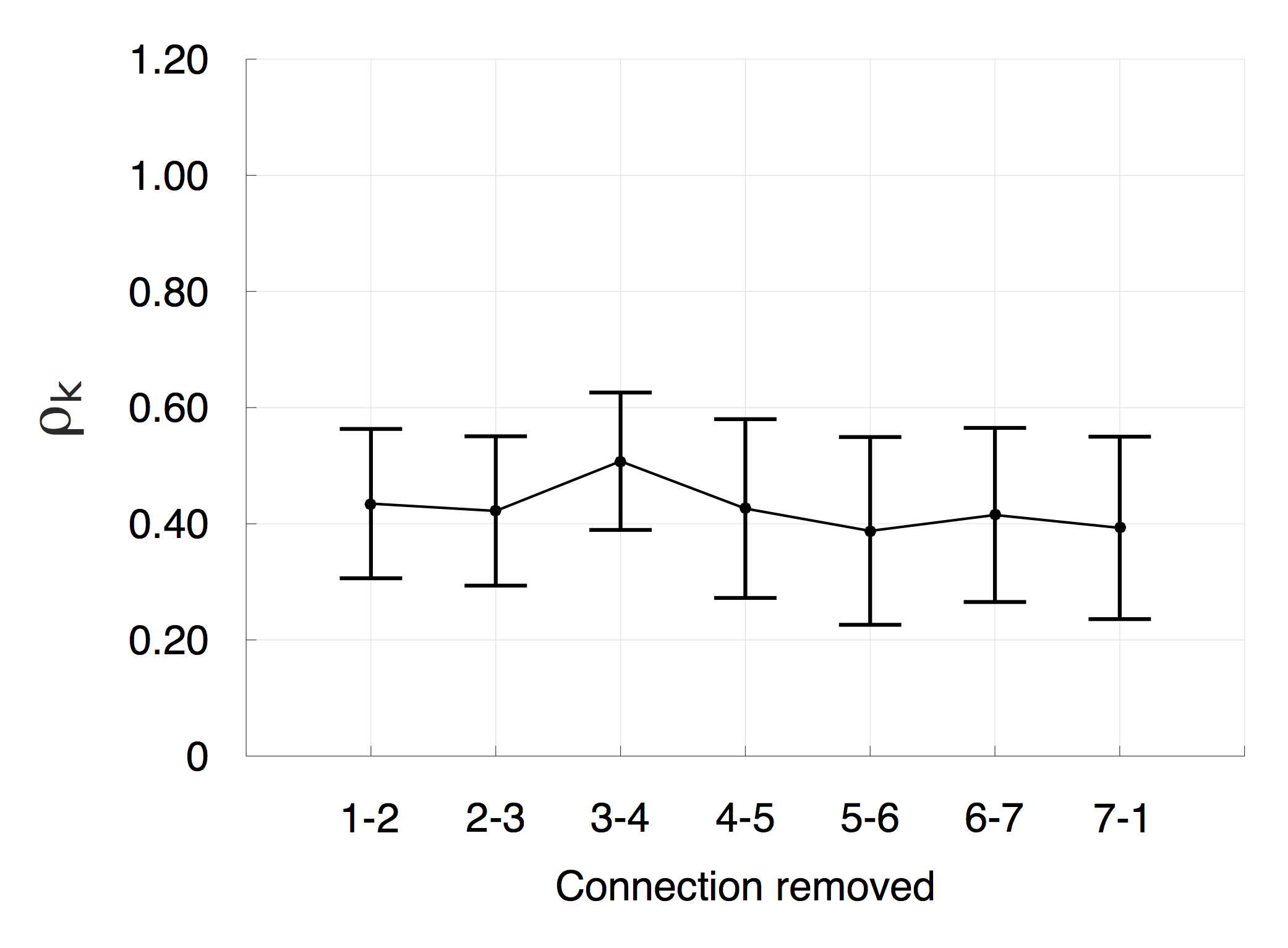}
 \caption{\textbf{Individual synchronisation indices $\rho_k$ as a function of the link getting removed in a Ring graph topology to form a Path Graph.} Mean values over the total number of nodes are represented by circles, and standard deviations by error bars.}
 \label{fig:S1_fig11}
\end{figure}

\subsection{Selection of the central node in a Star graph}

We finally considered a network of $N=7$ heterogeneous Kuramoto oscillators connected over a Star graph topology (Fig. 2d in the main text). Seven scenarios were considered, where each scenario differs from the others in the selection of the central node. Specifically, in the first scenario the central node was set to be node $1$, in the second scenario node $2$, up to the seventh scenario where the central node was set to be node $7$. In order to isolate the effects of such choice on the coordination levels $\rho_k$, the group dispersion was not varied over the different scenarios ($c_v=17\%$), and neither were the individual variabilities across all the nodes ($\sigma(\omega_k)=0.25 \ \forall k \in [1,N]$).  The other parameters were set as in the previous cases.

A One-way ANOVA using the Welch's test revealed that the differences in the coordination levels obtained for different choices of the central node in a Star graph topology are not statistically significant ($F(6,18.367)=0.463$, $p=0.827$, $\eta^2=0.070$). Mean values and standard deviations of the individual synchronisation index $\rho_k$ obtained in the seven different scenarios here considered are shown in Supplementary Fig. \ref{fig:S1_fig12}.

\begin{figure}[!ht]
 \centering
 \includegraphics[width=.75\textwidth]{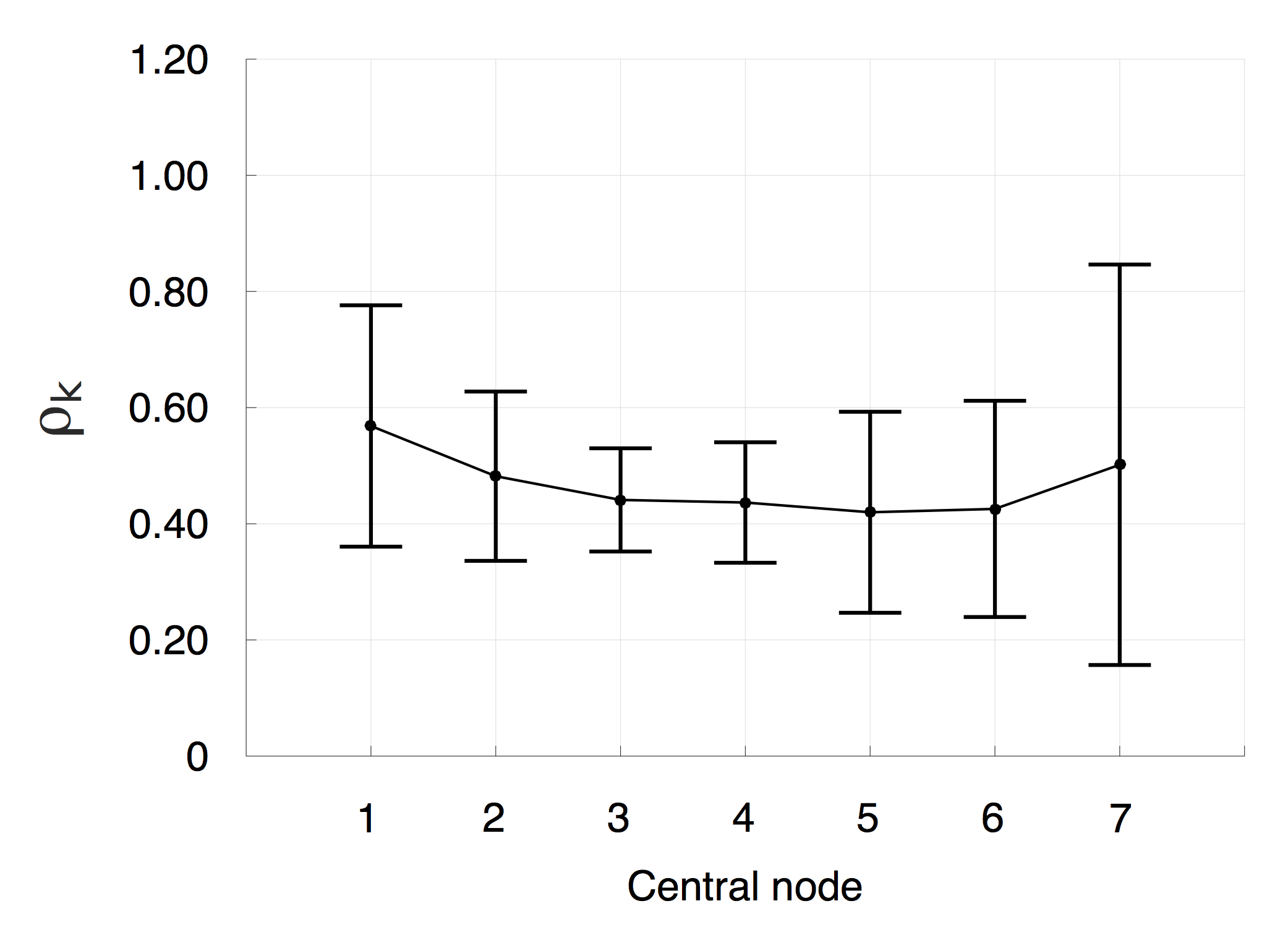}
 \caption{\textbf{Individual synchronisation indices $\rho_k$ as a function of the central node selection in a Star graph topology.} Mean values over the total number of nodes are represented by circles, and standard deviations by error bars.}
 \label{fig:S1_fig12}
\end{figure}



\begin{thebibliography}{9}

\bibitem{NRBVB00}
N{\'e}da, Z., Ravasz, E., Brechet, Y., Vicsek, T. \& Barab{\'a}si, A.-L. Self-organizing processes: The sound of many hands clapping. \emph{Nature} \textbf{403(6772),} 849-850 (2000).

\bibitem{MPGHT10}
Moussa{\"\i}d, M., Perozo, N., Garnier, S., Helbing, D. \& Theraulaz, G. The walking behaviour of pedestrian social groups and its impact on crowd dynamics. \emph{PLoS One} \textbf{5(4),} e10047 (2010).

\bibitem{RW14}
Rio, K. \& Warren, W. H. The visual coupling between neighbors in real and virtual crowds. \emph{Transport. Res. Proc.} \textbf{2,} 132-140 (2014).

\bibitem{DBLTCCAF12}
D'Ausilio, A. \emph{et al.} Leadership in orchestra emerges from the causal relationships of movement kinematics. \emph{PLoS One} \textbf{7(5),}  e35757 (2012).

\bibitem{BDGCF14}
Badino, L., D'Ausilio, A., Glowinski, D., Camurri, A. \& Fadiga, L. Sensorimotor communication in professional quartets. \emph{Neuropsychologia} \textbf{55,} 98-104 (2014).

\bibitem{DWA10}
Duch, J., Waitzman, J. S. \& Amaral, L. A. N. Quantifying the performance of individual players in a team activity. \emph{PLoS One} \textbf{5(6),} e10937 (2010).

\bibitem{S16}
Silva, P. \emph{et al.} Practice effects on intra-team synergies in football teams. \emph{Hum. Movement Sci.} \textbf{46,} 39-51 (2016).

\bibitem{L14}
Leonard, N. E. \emph{et al.} \emph{In the Dance Studio: An Art and Engineering Exploration of Human Flocking}. Springer, Controls and Art, 27-49 (2014).
  
\bibitem{DLTL16}
Dey, B., Lazier, R. J., Trueman, D. \& Leonard, N. E. Investigating group behavior in dance: an evolutionary dynamics approach. \emph{Proc. IEEE American Control Conf.} Boston, Massachusetts, USA. 6465-6470 (2016). 

\bibitem{EBQBM16}
Ellamil, M., Berson, J., Wong, J., Buckley, L. \& Margulies, D. S. One in the dance: Musical correlates of group synchrony in a real-world club environment. \emph{PLoS One} \textbf{11(10),} e0164783 (2016).

\bibitem{VRHKV04}
Van Baaren, R. B., Holland, R. W., Kawakami, K. \& Van Knippenberg, A. Mimicry and prosocial behavior. \emph{Psychol. Sci.} \textbf{15(1),} 71-74 (2004).

\bibitem{WH09}
Wiltermuth, S. S. \& Heath, C. Synchrony and cooperation. \emph{Psychol. Sci.} \textbf{20(1),} 1-5 (2009).

\bibitem{CKFL05}
Couzin, I. D., Krause, J., Franks, N. R. \& Levin, S. A. Effective leadership and decision-making in animal groups on the move. \emph{Nature} \textbf{433(7025),}  513-516 (2005).

\bibitem{NABV10}
Nagy, M., \'{A}kos, Z., Biro, D. \& Vicsek, T. Hierarchical group dynamics in pigeon flocks. \emph{Nature} \textbf{464(7290),}  890-893 (2010).

\bibitem{NVPMVB13}
Nagy, M. \emph{et al.} Context-dependent hierarchies in pigeons. \emph{Proc. Natl. Acad. Sci. U.S.A.} \textbf{110(32),}  13049-13054 (2013).

\bibitem{ZBPdB15}
Zienkiewicz, A., Barton, D. A. W., Porfiri, M. \& di Bernardo, M. Leadership emergence in a data-driven model of zebrafish shoals with speed modulation. \emph{Eur. Phys. J. Spec. Top.} \textbf{224(17-18),} 3343-3360 (2015).

\bibitem{OdGJLK08}
Oullier, O., De Guzman, G. C., Jantzen, K. J., Lagarde, J. \& Scott Kelso, J. A. Social coordination dynamics: Measuring human bonding. \emph{Soc. Neurosci.} \textbf{3(2),} 178-192 (2008).

\bibitem{ST94}
Schmidt, R. C. \& Turvey, M. T. Phase-entrainment dynamics of visually coupled rhythmic movements. \emph{Biol. Cybern.} \textbf{70(4),} 369-376 (1994).

\bibitem{VMLB11}
Varlet, M., Marin, L., Lagarde, J. \& Bardy, B. G. Social postural coordination. \emph{J. Exp. Psychol. Hum. Percept. Perform.} \textbf{37(2),} 473-483 (2011).

\bibitem{slowinski2015dynamic}
S\l{}owi\'{n}ski, P. \emph{et al.} Dynamic similarity promotes interpersonal coordination in joint action. \emph{J. R. Soc. Interface} \textbf{13(116),} 20151093 (2016)

\bibitem{NDA11}
Noy, L., Dekel, E. \& Alon, U. The mirror game as a paradigm for studying the dynamics of two people improvising motion together. \emph{Proc. Natl. Acad. Sci. USA} \textbf{108(52),} 20947-20952 (2011).

\bibitem{ZATdMSMC}
Zhai, C., Alderisio, F., Tsaneva-Atanasova, K. \& di Bernardo, M. A novel cognitive architecture for a human-like virtual player in the mirror game. \emph{Proc. IEEE Conf. Syst., Man, Cybern.} San Diego, California, USA. 754-759 (2014).

\bibitem{ZASTdB16}
Zhai, C., Alderisio, F., S\l{}owi\'{n}ski, P., Tsaneva-Atanasova, K. \& di Bernardo, M. Design of a virtual player for joint improvisation with humans in the mirror game. \emph{PLoS One} \textbf{11(4),} e0154361 (2016).

\bibitem{FR10}
Frank, T. D. \& Richardson, M. J. On a test statistic for the Kuramoto order parameter of synchronization: an illustration for group synchronization during rocking chairs. \emph{Phys. D} \textbf{239(23),}  2084-2092 (2010).

\bibitem{RGFGM12}
Richardson, M. J., Garcia, R. L., Frank, T. D., Gergor, M. \& Marsh, K. L. Measuring group synchrony: a cluster-phase method for analyzing multivariate movement time-series. \emph{Front. Physiol.} \textbf{3} (2012).

\bibitem{ABdB15}
Alderisio, F., Bardy, B. G. \& di Bernardo, M. Entrainment and synchronization in networks of Rayleigh--van der Pol oscillators with diffusive and Haken--Kelso--Bunz couplings. \emph{Biol. Cybern.} \textbf{110(2),} 151-169, DOI: 10.1007/s00422-016-0685-7 (2016).

\bibitem{IR15} 
Iqbal, T. \& Riek, L. A method for automatic detection of psychomotor entrainment. \emph{IEEE T. Affect. Comput.} \textbf{7(1),} 3-16 (2016)

\bibitem{HT09} 
Himberg, T. \& Thompson, M. Group synchronization of coordinated movements in a cross-cultural choir workshop. \emph{7th Triennial Conference of European Society for the Cognitive Sciences of Music} 175-180 (2009).

\bibitem{CBVB14} 
Codrons, E., Bernardi, N. F., Vandoni, M. \& Bernardi, L. Spontaneous group synchronization of movements and respiratory rhythms. \emph{PLoS One} \textbf{9(9),} e107538 (2014).

\bibitem{WW95} 
Wing, A. M. \& Woodburn, C. The coordination and consistency of rowers in a racing eight. \emph{J. Sports Sci.} \textbf{13(3),} 187-197 (1995).

\bibitem{YY11}
Yokoyama, K. \& Yamamoto, Y. Three people can synchronize as coupled oscillators during sports activities. \emph{PLoS Comput. Biol.} \textbf{7(10),} e1002181 (2011).

\bibitem{healey2005inter}
Healey, P. G. T., Leach, J. \& Bryan-Kinns N. Inter-play: Understanding group music improvisation as a form of everyday interaction. \emph{Proceedings of Less is More--Simple Computing in an Age of Complexity} (2005).

\bibitem{kauffeld2009complaint}
Kauffeld, S. \& Meyers, R. A. Complaint and solution-oriented circles: Interaction patterns in work group discussions. \emph{Eur. J. Work Organ. Psy.} \textbf{18(3),} 267-294 (2009).

\bibitem{passos2011networks}
Passos, P. \emph{et al.} Networks as a novel tool for studying team ball sports as complex social systems. \emph{J. Sci. Med. Sport} \textbf{14(2),} 170-176 (2011).

\bibitem{duarte2012intra}
Duarte, R. \emph{et al.} Intra-and inter-group coordination patterns reveal collective behaviors of football players near the scoring zone. \emph{Hum. Mov. Sci.} \textbf{31(6),} 1639-1651 (2012).

\bibitem{RTAR13}
Duarte, R., Travassos, B., Ara{\'u}jo, D. \& Richardson, M. J. \emph{The influence of manipulating the defensive playing method on collective synchrony of football teams}. Performance Analysis of Sport IX. Routledge, Taylor \& Francis Group London (2013).

\bibitem{K84}
Kuramoto, Y. \emph{Chemical Oscillations, Waves and Turbulence}. Springer, Heidelberg (1984).

\bibitem{strogatz2000kuramoto}
Strogatz, S. H. From Kuramoto to Crawford: exploring the onset of synchronization in populations of coupled oscillators. \emph{Physica D} \textbf{143(1)} 1-20 (2000).

\bibitem{BL95}
Baumeister, R. F. \& Leary, M. R. The need to belong: desire for interpersonal attachments as a fundamental human motivation. \emph{Psychol. Bull.} \textbf{117(3),} 497-529 (1995).

\bibitem{MFH10}
M{\"a}s, M., Flache, A. \& Helbing, D. Individualization as driving force of clustering phenomena in humans. \emph{PLoS Comput. Biol.} \textbf{6(10),} e1000959 (2010).

\bibitem{SFV13}
Stark, T. H., Flache, A. \& Veenstra, R. Generalization of positive and negative attitudes toward individuals to outgroup attitudes. \emph{Pers. Soc. Psychol. B.} \textbf{39(5),} 608-622 (2013).

\bibitem{VDBCF16}
Volpe, G., D'Ausilio, A., Badino, L., Camurri, A. \& Fadiga, L. Measuring social interaction in music ensembles. \emph{Phil. Trans. R. Soc. B.} \textbf{371,} 20150377 (2016).

\bibitem{STdB15}
Scafuti, F., Takaaki, A. \& di Bernardo, M. Heterogeneity induces emergent functional networks for synchronization. \emph{Phys. Rev. E} \textbf{91(6),} 062913 (2015).

\bibitem{AB02}
Albert, R. \&  Barab{\'a}si, A.-L. Statistical mechanics of complex networks. \emph{Rev. Mod. Phys.} \textbf{74(1)}, 47-97 (2002).

\bibitem{NBW11}
Newman, M., Barab{\'a}si, A.-L \& Watts, D. J. The structure and dynamics of networks. \emph{Princeton University Press} (2011).

\bibitem{ADKMZ08}
Arenas, A., D{\'\i}az-Guilera, A., Kurths, J., Moreno, Y. \& Zhou, C. Synchronization in complex networks. \emph{Phys. Rep.} \textbf{469(3),} 93-153 (2008).

\bibitem{AGGLB11}
Assenza, S., Guti{\'e}rrez, R., G{\'o}mez-Garde{\~n}es, J., Latora, V. \& Boccaletti, S. Emergence of structural patterns out of synchronization in networks with competitive interactions. \emph{Sci. Rep.} \textbf{1} (2011).

\bibitem{MP04}
Moreno, Y. \& Pacheco, A. F. Synchronization of Kuramoto oscillators in scale-free networks. \emph{Europhys. Lett.} \textbf{68(4),} 603 (2004).

\bibitem{SRP06}
Santos, F. C., Rodrigues, J. F. \& Pacheco, J. M. Graph topology plays a determinant role in the evolution of cooperation. \emph{Proc. R. Soc. Lond. [Biol]} \textbf{273(1582),} 51-55 (2006).

\bibitem{AC16}
Antonioni, A. \& Cardillo, A. Coevolution of synchronization and cooperation in costly networked interactions. \emph{Phys. Rev. Lett.} \textbf{118(23),} 238301 (2017).

\bibitem{ZWP16}
Zamm, A., Wellman, C. \& Palmer, C. Endogenous rhythms influence interpersonal synchrony. \emph{J. Exp. Psychol. Hum. Percept. Perform.} 1-6 (2016).

\bibitem{ALFdB16}
Alderisio, F., Lombardi, M., Fiore, G., \& di Bernardo, M. Study of movement coordination in human ensembles via a novel computer-based set-up. arXiv preprint arXiv:1608.04652 (2016).

\bibitem{BSMGC16}
Boucenna, S., Cohen, D., Meltzoff, A. N., Gaussier, P. \& Chetouani, M. Robots learn to recognize individuals from imitative encounters with people and avatars. \emph{Sci. Rep.} \textbf{6} (2016).

\bibitem{IRR16}
Iqbal, T., Rack, S. \& Riek, L. D. Movement coordination in human-robot teams: a dynamical systems  approach. arXiv preprint arXiv:1605.01459 (2016).

%

\bibitem{AAZdB16}
Alderisio, F., Antonacci, D., Zhai, C., \& di Bernardo, M. Comparing different control approaches to implement a human-like virtual player in the mirror game. \emph{Proc. European Control Conference}, Aalborg, Denmark. 216-221 (2016)

\bibitem{ZASTdB17}
Zhai, C., Alderisio, F., S\l{}owi\'{n}ski, P., Tsaneva-Atanasova, K., \& di Bernardo, M. Design and validation of a virtual player for studying interpersonal coordination in the mirror game. \emph{IEEE Trans. Cybern.}, 10.1109/TCYB.2017.2671456 (2017)

\bibitem{MRS09}
Marsh, K. L., Richardson, M. J. \& Schmidt, R. C. Social connection through joint action and interpersonal coordination. \emph{Top. Cogn. Sci.} \textbf{1(2),} 320-339 (2009).

\bibitem{IHV02}
Farkas, I., Helbing, D. \& Vicsek, T. Social behaviour: Mexican waves in an excitable medium. \emph{Nature} \textbf{419(6903),} 131-132 (2002).

\bibitem{KCRPM08}
Kralemann, B., Cimponeriu, L., Rosenblum, M., Pikovsky, A. \& Mrowka, R. Phase dynamics of coupled oscillators reconstructed from data. \emph{Phys. Rev. E} \textbf{77(6),}  066205 (2008).

\bibitem{dLdBL15}
DeLellis, P., di Bernardo, M. \& Liuzza, D. Convergence and synchronization in heterogeneous networks of smooth and piecewise smooth systems. \emph{Automatica} \textbf{56,} 1-11 (2015).

\end{thebibliography}
\end{document}